\documentclass{amsart}
\usepackage[frenchb, english]{babel}

\newtheorem{theoreme}{Th\'eoreme}
\newtheorem*{theoreme*}{Th\'{e}or\`{e}me}
\newtheorem{lemme}{Lemme}
\newtheorem{corollaire}{Corollaire}
\newtheorem*{corollaire*}{Corollaire}
\newtheorem{proposition}{Proposition}
\newtheorem*{proposition*}{Proposition}
\theoremstyle{definition}
\newtheorem*{definition*}{D\'efinition}
\newtheorem{definition}{D\'efinition}

\theoremstyle{remark}
\newtheorem*{remarque}{Remarque}
\newtheorem*{notation}{Notation}
\newtheorem*{notations}{Notations}
\theoremstyle{plain}

\newcommand{\Z}{{\mathbb Z}}
\newcommand{\T}{{\mathbb T}}
\newcommand{\C}{{\mathbb C}}
\newcommand{\CZ}{{\mathcal Z}}
\newcommand{\bn}{\mathbf{n}}
\newcommand{\bx}{\mathbf{x}}
\newcommand{\by}{\mathbf{y}}
\newcommand{\bj}{\mathbf{j}}
\newcommand{\bz}{\mathbf{z}}
\newcommand{\bu}{\mathbf{u}}
\newcommand{\bv}{\mathbf{v}}
\newcommand{\bw}{\mathbf{w}}
\newcommand{\ubx}{\mathbf{\underline x}}

\newcommand{\ubz}{\mathbf{\underline z}}
\DeclareMathOperator{\CG}{\textsf{\upshape G}}
\DeclareMathOperator{\CP}{\textsf{\upshape P}}
\DeclareMathOperator{\CQ}{\textsf{\upshape Q}}

\newcommand{\inv}{^{-1}}
\DeclareMathOperator{\id}{\bf Id}
\newcommand{\equivq}{\underset{\CQ}{\equiv}}
\newcommand{\type}[1]{^{[#1]}}
\newcommand{\br}[1]{\left[#1\right]}
\newcommand{\CGtop}{\CG_{\text{\small\upshape top}}}
\newcommand{\Gammatop}{\Gamma_{\text{\small\upshape top}}}
\DeclareMathOperator{\RP}{\bf RP}
\DeclareMathOperator{\RPD}{\bf RP^2}
\DeclareMathOperator{\RPDS}{\bf RP^2_{\!s}}

\begin{document}
\selectlanguage{frenchb}
\def\and{et}

\title[Nilsyst\`emes d'ordre deux]
{Nilsyst\`emes d'ordre deux et  parall\'el\'epip\`edes}
\author{Bernard Host \& Alejandro Maass}
\address{B. Host : \'Equipe d'analyse et de math\'ematiques appliq\'ees, Universit\'e 
de Marne la Vall\'ee\\
5 Bd. Desacrtes, Champs sur Marne\\
77454 Marne la Vall\'ee cedex -- France}
\email{bernard.host\@ univ-mlv.fr}

\address{A. Maass : Departamento de Ingenier\'{\i}a
Matem\'atica, Universidad de Chile
\& Centro de Modelamiento Ma\-te\-m\'a\-ti\-co,
UMI 2071 UCHILE-CNRS, Casilla 170/3 correo 3,
Santiago, Chile.}
\email{amaass@dim.uchile.cl} 
\def\datename{Date :}
\date{2 Ao\^ut 2006}
\maketitle

\begin{abstract}
En topologie dynamique, une famille classique de syst\`emes est celle 
form\'ee par les rotations minimales. La classe des 
nilsyst\`emes et de leurs limites 
projectives en est une extension naturelle. 
L'\'etude de ces syst\`emes est ancienne mais conna\^\i t actuellement un 
renouveau \`a cause de ses applications, \`a la fois \`a la th\'eorie ergodique
et en th\'eorie additive des nombres.

Les rotations minimales sont caract\'eris\'ees par le fait que la relation 
de proximalit\'e r\'egionale est l'\'egalit\'e. Nous introduisons une nouvelle 
relation, celle de \emph{bi-proximalit\'e r\'egionale}, et montrons qu'elle 
caract\'erise les limites projectives de nilsyst\`emes d'ordre deux.

Les rotations minimales sont li\'ees aux suites presque p\'eriodiques et de 
m\^eme les nilsyt\`emes correspondent aux \emph{nilsuites}. Ces suites 
introduites en th\'eorie ergodique sont intervenues depuis dans 
certaines questions de th\'eorie des nombres.
De notre caract\'erisation des nilsyst\`emes d'ordre 
deux nous d\'eduisons une caract\'erisation des nilsuites d'ordre deux.

Les d\'emonstrations s'appuient d'une mani\`ere essentielle sur l'\'etude 
des 
\guillemotleft~structures de parall\'el\'epip\`edes~\guillemotright\ 
d\'evelopp\'ee  par B. Kra 
et le 
premier auteur.
\end{abstract}

\selectlanguage{english}
\begin{abstract}
A classic family in topological dynamics is that of minimal rotations.
One natural extension of this family is the class of nilsystems and
their inverse limits. 
These systems have arisen in recent applications
in ergodic theory and in additive combinatorics, renewing interest in
studying these classical objects. 

Minimal rotations can be
characterized via the regionally proximal relation. 
We introduce a new
relation, the bi-regionally proximal relation, and show that it
characterizes inverse limits of two step nilsystems. 

Minimal rotations are linked to almost periodic sequences, 
and more generally nilsystems correspond to nilsequences.
Theses sequences were introduced in ergodic theory and have since be 
used in some questions of Numer Theory. 
Using our characterization of two step nilsystems  we
deduce a characterization of two step nilsequences.

The proofs rely on in an essential way  the study of ``parallelepiped 
structures'' developed  by B. Kra and the first author.
\end{abstract}

\selectlanguage{frenchb}

\section{Introduction}
Un \emph{syst\`eme dynamique topologique} (abr\'eg\'e en \guillemotleft~syst\`eme~\guillemotright)
$(X,T)$ est un espace compact m\'etrisable $X$ muni d'un hom\'eomorphisme
$T\colon X\to X$. 
Ce syst\`eme est \emph{transitif} s'il existe au moins 
un point $x\in X$ dont l'\emph{orbite}  $\{T^nx\;;\;n\in\Z)$ est dense 
dans $X$ ; 
il est \emph{minimal} si l'orbite de tout point est dense 
dans $X$. 
Cette propri\'et\'e est \'equivalente \`a la condition que les seuls ferm\'es de
$X$ invariants par $T$ sont $X$ et l'ensemble vide.

\subsection{Syst\`emes \'equicontinus et nilsyst\`emes}

Une famille classique  de syst\`emes minimaux est celle form\'ee par les 
\emph{rotations minimales} (ces syt\`emes sont appel\'es 
\emph{syst\`emes de Kronecker} dans~\cite{Fu}). Dans ce cas, $X$ est un groupe ab\'elien 
compact et la transformation $T$ a la forme $x\mapsto \alpha x$, o\`u 
$\alpha\in X$ est une constante telle que 
$\{\alpha^n\;;\;n\in\Z\}$ 
soit dense dans $X$. 

Un syst\`eme minimal est \emph{\'equicontinu} si et 
seulement si il est isomorphe \`a une rotation minimale. 
Ces syst\`emes peuvent \^etre caract\'eris\'es au moyen de la
\emph{relation de proximalit\'e r\'egionale} : pour chaque syst\`eme 
minimal $(X,T)$, cette relation est une relation d'\'equivalence sur $X$ 
et le syst\`eme est \'equicontinu si et seulement si cette relation est 
l'\'egalit\'e~\cite{Aus}. La d\'efinition et les propri\'et\'es de cette 
relation sont rappel\'ees dans la section~\ref{subsec:RP}.

Nous nous int\'eressons ici \`a une famille plus g\'en\'erale de syst\`emes, 
les \emph{nilsyst\`emes}. 
L'\'etude de ces syst\`emes est ancienne mais conna\^\i t actuellement un 
renouveau \`a cause de ses application, \`a la fois en th\'eorie ergodique
 (voir par exemple~\cite{CL2}, \cite{CL}, \cite{HK2}, \cite{BHK}) 
et en th\'eorie additive des nombres (voir~\cite{GT}, \cite{GT2}).
Nous en donnons ici la d\'efinition, les propri\'et\'es utiles sont 
rappel\'ees dans la section~\ref{subsec:rappels-nil}.

\begin{definition}
\label{def:nilsysteme}
Soient $k\geq 1$ un entier, $G$ un groupe de Lie nilpotent d'ordre $k$
 et $\Gamma$ un sous-groupe discret cocompact de $G$.

La vari\'et\'e $X=G/\Gamma$ est appel\'ee une 
\emph{nilvari\'et\'e d'ordre $k$}.

Fixons $t\in G$ et soit $T\colon X\to X$ l'application 
$x\mapsto t\cdot x$, 
o\`u $(g,x)\mapsto g\cdot x$ est l'action \`a gauche de $G$ sur
$X$. 
Alors $(X,T)$ est appel\'e un \emph{nilsyst\`eme d'ordre $k$}.
\end{definition}
On rappelle qu'un nilsyst\`eme est minimal d\`es qu'il est transitif.

Un nilsyst\`eme d'ordre $1$ est une rotation sur un groupe 
ab\'elien compact. Nous ne consid\'erons ici que les nilsyst\`emes d'ordre $2$, la g\'en\'eralisation de nos r\'esultats \`a l'ordre sup\'erieur 
restant une question ouverte.

Contrairement \`a la famille des rotations minimales, la famille des 
nilsyst\`{e}mes minimaux d'ordre deux n'est pas stable par limite projective
(voir~\cite{Ru}).
Il est donc naturel de consid\'erer la classe plus large des syst\`emes 
qui sont des \emph{limites projectives de nilsyst\`emes minimaux 
d'ordre deux}.

Le r\'esultat principal de cet article (th\'eor\`eme~\ref{th:main}) est
 la caract\'erisation des limites projectives de
nilsyst\`emes minimaux d'ordre deux au moyen d'une nouvelle relation, 
la \emph{relation de bi-proximalit\'e r\'egionale} introduite dans la 
section~\ref{subsec:RP} : Cette relation est l'\'egalit\'e si et 
seulement si le syst\`eme est (isomorphe \`a) une limite projective de
 nilsyst\`eme d'ordre deux.

\subsection{Application aux nilsuites}

La notion de  \emph{nilsuite} (section~\ref{subsec:def-nilsuites})
est une g\'en\'e\-ra\-lisation naturelle de 
celle de suite presque-p\'eriodique, obtenue en rempla\c cant les 
groupes ab\'eliens par des groupes nilpotents. Ces suites ont \'et\'e 
d'abord introduites en th\'eorie ergodique pour dans l'\'etude des 
corr\'elations multiples (\cite{BHK}) et elles ont \'et\'e ensuites 
utilis\'ees en th\'eorie des nombres pour l'\'etude des configurations 
apparaissant dans les nombres premiers (\cite{GT,GT2}).

Dans la section~\ref{sec:nilsuites} nous utilisons notre r\'{e}sultat 
principal pour 
donner une caract\'erisation des nilsuites d'ordre deux en termes
 de r\'egularit\'e arithm\'etique. Nous pensons que cette 
caract\'erisation pourrrait avoir des applications en dehors de la 
dynamique. 
Au passage nous donnons aussi une 
 caract\'erisation des suites presque p\'eriodiques, facile \`a 
 montrer directement mais 
 apparemment nouvelle.

\subsection{Facteur \'equicontinu maximal et nilfacteur d'ordre deux
maximal}

En topologie dynamique, les quotients sont appel\'es des 
facteurs. 
Nous en rappelons la d\'efinition.

\begin{definition}
Soit  $(X,T)$ un syst\`eme.
Un \emph{facteur} de ce 
syst\`eme est la donn\'ee d'un syst\`eme $(Y,S)$ et d'une application
 $p\colon X\to Y$, continue, surjective et v\'erifiant
 $S\circ p=p\circ T$. 
\end{definition}

Par abus de langage, on dit bri\`evement que $Y$ est un facteur de $X$ 
et l'application $p$ est appel\'ee l'\emph{application facteur}.

Remarquons que tout facteur d'un syst\`eme minimal  est minimal et que 
tout facteur d'un syst\`eme minimal \'equicontinu est \'equicontinu.

Soit d\'esormais $(X,T)$ un syst\`eme minimal. 
Les facteurs de $X$ sont (partiellement) ordonn\'es de 
fa\c con naturelle : Si $Y$ est un facteur de $X$ et si $Z$ est un facteur de $Y$ alors 
$Z$ est un facteur de $X$ et on dit alors que la facteur $Y$ est 
\emph{plus grand} ou \emph{au dessus} du facteur $Z$ de $X$.

La famille des facteurs de $X$ qui 
sont \'equicontinus est projective. Comme une limite projective
de syst\`emes \'equicontinus est un syst\`eme \'equicontinu, cette famille a 
un plus grand \'el\'ement, le \emph{facteur \'equicontinu maximal}. 
Cette notion est l'analogue en topologie 
dynamique ce celle de \emph{facteur de Kronecker} en th\'eorie 
ergodique.
\begin{notation}
Nous notons d\'esormais $Z$ ou $Z_X$ le facteur \'equicontinu 
maximal de $(X,T)$ et $\pi\colon X\to Z$ ou $\pi_X\colon X\to Z_X$ l'application facteur.
\end{notation}

Le facteur \'equicontinu maximal est d\'etermin\'e par la relation de 
proximalit\'e r\'egionale : Cette relation est une relation 
d'\'equivalence ferm\'ee et invariante et $Z$ est le quotient de $X$ par 
cette relation (\cite{Aus}). 

Consid\'erons maintenant la famille des facteurs de $X$ qui sont des 
nilsyst\`emes d'ordre $2$. 
Cette famille est projective 
(proposition~\ref{prop:nilprojectif}), mais n'est pas stable par 
limite projective. 

 \begin{definition}
Le \emph{nilfacteur d'ordre $2$ maximal} de $X$ est le plus grand 
facteur de $X$ qui est (isomorphe \`a) une 
limite projective de nilsyst\`emes d'ordre $2$.
\end{definition}
\begin{notation}
Nous notons d\'esormais $Z_2$ ou $Z_{X,2}$ le nilfacteur d'ordre $2$
maximal de $(X,T)$ et $\pi_2$ ou $\pi_{X,2}$ l'application quotient.
\end{notation}

Ainsi, $Z_2$ est caract\'eris\'e par les deux propri\'et\'es :
\begin{enumerate}
\item
$Z_2$ est un facteur de $X$ qui est isomorphe \`a une limite 
projective de nilsyst\`emes d'ordre $2$ et
\item
Tout facteur de $X$ qui est isomorphe \`a un nilsyst\`eme d'ordre $2$  est un facteur de $Z_2$.
\end{enumerate}

La notion de nilfacteur d'ordre $2$
maximal est l'analogue en topologie  dynamique de l'\emph{alg\`ebre de 
Conze-Lesigne} construite dans~\cite{CL2} et \cite{CL} puis 
dans~\cite{HK2} o\`u elle est not\'ee $Z_2$.

Nous ne savons pas si en g\'en\'eral le nilfacteur d'ordre $2$ 
maximal est le quotient de $X$ par la relation de 
bi-proximalit\'e r\'egionale, 
ni m\^eme si cette relation est d'\'equivalence. 
Cependant ces propri\'et\'es sont vraies (th\'eor\`eme~\ref{th:main-distal}) 
pour une classe 
importante de syst\`emes, les \emph{syst\`emes distaux} dont la 
d\'efinition 
est rappel\'ee dans la section~\ref{subsec:def-distal}.

\subsection{Parall\'elogrammes et parall\'el\'epip\`edes dynamiques}

Soit $(X,T)$ un sys\-t\`eme minimal. 
Dans la section~\ref{subsec:def-para} nous construisons un 
sous-ensemble $\CP$ de $X^4$ et un sous-ensemble $\CQ$ de $X^8$, 
appel\'es respectivement \emph{l'ensemble des parall\'elogrammes} et
\emph{l'ensemble des parall\'el\'epip\`edes} de $X$. Ce sous-ensembles 
v\'erifient une partie importante des propri\'et\'es stipul\'{e}es dans 
les d\'efinitions de~\cite{HK} d'une structure de parall\'elogrammes et  d'une 
structure de parall\'el\'epip\`edes, mais apparemment pas toutes.
 La principale propri\'et\'e manquante est 
la \guillemotleft~transitivit\'e~\guillemotright : Les huit points obtenus en accolant deux 
parall\'{e}l\'{e}pip\`{e}des le long d'une face commune ne semblent pas en g\'en\'eral 
former un parall\'el\'epip\`ede.

Cependant nous montrons (section~\ref{sec:distaux}) que dans le cas d'un syst\`eme distal, $(\CP,\CQ)$ est une 
structure de parall\'el\'epip\`edes au sens de~\cite{HK}. Cela permet 
d'utiliser dans ce cas tous les r\'esultats de cet article. 
Pour traiter le cas g\'en\'eral il suffit alors de remarquer 
(proposition~\ref{prop:distalRP}) 
que si la 
relation de bi-proximalit\'e r\'egionale est triviale alors le syst\`eme est 
distal.

\subsection{}

Les d\'emonstrations s'inspirent de la
construction  de 
l'alg\`ebre de Conze-Lesigne dans~\cite{HK2}.
 Cependant les m\'ethodes sont tr\`es 
diff\'erentes, aucun des outils de base de la construction ergodique
(esp\'erance conditionnelle, cocycles, \'equations fonctionnelles, 
seminormes, \dots) n'ayant d'analogue en topologie dynamique. 

Les outils \guillemotleft~alg\'ebriques~\guillemotright\ de~\cite{HK} remplacent ici ces 
notions ergodiques, malgr\'{e} la difficult\'e caus\'ee par le fait que $\CQ$ 
n'est pas une vraie structure de parall\'el\'epip\`edes. 

Il est naturel d'essayer d'\'etendre nos r\'esultats aux ordres 
sup\'erieurs : pour tout syst\`eme minimal $(X,T)$ et pour chaque $k\geq 2$ 
on peut d\'efinir une  relation de r\'egionalement proximale d'ordre $k$ 
sur $X$ et un nilfacteur d'ordre $k$ maximal de $X$. On peut 
conjecturer qu'un syst\`eme est une limite projective de 
nilsyst\`emes d'ordre 
 $k$ si et seulement si la relation d'ordre $k$ est l'\'egalit\'e. 
Pour attaquer ce probl\`eme il faudrait d'abord g\'en\'eraliser les 
r\'esultats alg\'ebriques de~\cite{HK} \`a l'ordre sup\'erieur. Une des 
principales difficult\'es techniques est la m\^eme que dans le cadre de la 
th\'eorie 
ergodique : Les syst\`emes \`a consid\'{e}rer \`a un niveau $k$ ne sont pas 
seulement les nilsyst\`emes d'ordre $k$ mais les limites projectives 
de tels syst\`emes, ce qui rend les r\'ecurrences beaucoup plus 
difficiles \`a mettre en \oe uvre.

Par ailleurs, on  peut se demander si la m\'ethode de cet article 
ne pourrait pas de produire une d\'emonstration alternative et plus 
simple des r\'esultats de~\cite{HK}.

\section{R\'esultats}
Dans  ce qui suit, $(X,T)$ est un syst\`eme. $X$ est muni 
d'une distance not\'ee $d$ ; pour $x\in X$ et $\epsilon>0$, 
$B(x,\epsilon)$ d\'esigne la boule ouverte de centre $x$ et de rayon $\epsilon$. 

\subsection{La relation de double proximalit\'e r\'egionale}
\label{subsec:RP}
Rappelons la d\'efinition de la relation de proximalit\'e r\'egionale, 
\'etudi\'ee par Auslander dans~\cite{Aus}.
\begin{definition}
Deux points $x,y\in X$ sont \emph{r\'egionalement proximaux} si pour 
tout $\epsilon>0$ il existe $x',y'\in X$ et $n\in\Z$ avec
$$
 d(x',x)<\epsilon\ ;\ d(y',y)<\epsilon\text{ et 
}d(T^nx',T^ny')<\epsilon\ .
$$
On note $\RP$ cette relation et \'egalement son graphe, c'est \`a dire l'ensemble 
des couples $(x,y)\in X^2$ tels que $x$ et $y$ soient r\'egionalement 
proximaux. En cas d'ambigu•t\'e on note $\RP(X)$ au lieu de $\RP$.
\end{definition}
On a :
\begin{theoreme}
[\cite{Aus}, chapitre 9]
\label{th:Aus}
Supposons que le syst\`eme $(X,T)$ est minimal. Alors
la relation de proximalit\'e r\'egionale est une relation d'\'equivalence 
ferm\'ee sur $X$ et 
le quotient de $X$ par cette relation est le facteur \'equicontinu 
maximal de $X$.
\end{theoreme}
En particulier, un syst\`eme minimal est \'equicontinu si et 
seulement si la relation de proximalit\'e r\'egionale est l'\'egalit\'e.

Nous d\'efinissons :
\begin{definition} \label{def:RPdeux}
Deux points $x,y\in X$ sont \emph{doublement r\'egionalement proximaux} si pour 
tout $\epsilon>0$ il existe $x',y'\in X$ et $m,n\in\Z$ avec
\begin{gather*}
d(x',x)<\epsilon\ ;\ d(y',y)<\epsilon\ ;\\
d(T^mx',T^my')<\epsilon\ ;\ d(T^nx',T^ny')<\epsilon
\text{ et }d(T^{m+n}x',T^{m+n}y')<\epsilon .
\end{gather*}
On note $\RPD$ (ou $\RPD(X)$) cette relation et \'egalement son graphe. 
\end{definition}
On remarque imm\'{e}diatement que cette relation est ferm\'ee
et plus fine que la relation de proximalit\'e r\'egionale (c'est \`a dire 
que
 $\RPD$ est inclus dans $\RP$).

Notre r\'esultat principal est :
\begin{theoreme}
\label{th:main}
Supposons que le syst\`eme $(X,T)$ est  transitif. Alors
 $(X,T)$ est une limite projective de 
nilsyst\`{e}mes 
d'ordre $2$ si et seulement si la relation de double proximalit\'e r\'egionale
 est l'\'egalit\'e.
\end{theoreme}
La plus grande partie de cet article est consacr\'ee \`a la d\'emonstration 
de ce  th\'eor\`eme, qui s'ach\`eve dans la section~\ref{sec:dyna}.

Il est commode d'introduire une variante dissym\'etrique de la relation $\RPD$.
\begin{definition}
Soient $x,y\in X$. On dit que ces points sont \emph{fortement doublement 
r\'egionalement proximaux} si pour tout $\epsilon>0$ il existe $x',y'\in 
X$ et $m,n\in \Z$ avec
\begin{gather*}
d(x',x)<\epsilon\ ;\ d(y',y)<\epsilon\ ;\\
d(T^mx',y)<\epsilon\ ;\ d(T^nx',y)<\epsilon\ ;\ 
d(T^{m+n}x',y)<\epsilon\\
d(T^my',y)<\epsilon\ ;\ d(T^ny',y)<\epsilon\text{ et }
d(T^{m+n}y',y)<\epsilon\ .
\end{gather*}
On note $\RPDS$ (ou $\RPDS(X)$) cette relation et \'egalement son graphe. 
\end{definition}
La relation $\RPDS$ est clairement plus fine que la relation $\RPD$. 
Cependant, on peut remplacer la relation $\RPD$ 
par la relation $\RPDS$
dans le th\'eor\`eme~\ref{th:main}. En effet, 
dans la section~\ref{subsec:RPDRPDS} nous montrons  :
\begin{proposition}
\label{prop:RPDRPDS}
Soit $(X,T)$ un syst\`eme transitif. Alors la relation $\RPDS$ est 
l'\'egalit\'e si et seulement si la relation $\RPD$ est l'\'egalit\'e.
\end{proposition}

\begin{theoreme}
\label{th:main-distal}
Supposons que le syst\`eme $(X,T)$ est distal et transitif (et donc minimal, 
voir la proposition~\ref{prop:aus-distaux}).
 Alors la relation 
de bi-proximalit\'e r\'egionale sur $X$ est une relation d'\'equivalence 
ferm\'ee et le quotient de $X$ par cette relation est le nilfacteur 
d'ordre deux maximal de $X$.

De plus, la relation $\RPDS$ co\"\i ncide avec la relation $\RPD$.
\end{theoreme}
Une partie de ce th\'eor\`eme est montr\'ee dans la section~\ref{subsec:RPDRPDS},
 le reste 
\'etant d\'eduit du th\'eor\`eme~\ref{th:main} dans la section~\ref{subsec:proof-main-distal}.

\subsection{Parall\'elogrammes et parall\'el\'epip\`edes dynamiques}
\label{subsec:def-para}

Les relations de proximalit\'e simple et double sont li\'ees \`a certains 
sous-ensembles de $X^4$ et de $X^8$ que nous introduisons maintenant.
\begin{definition}
\label{def:CPCQ}
L'\emph{ensemble des parall\'elogrammes} $\CP$ de $X$ est
l'adh\'erence dans $X^4$ de l'ensemble
$$
 \bigl\{(x,T^mx,T^nx,T^{m+n}x)\ ;\ x\in X,\ m,n\in\Z\bigr\}\ .
$$
L'\emph{ensemble des parall\'el\'epip\`edes} $\CQ$ de $X$ est
l'adh\'erence dans $X^8$ de l'ensemble
$$
 \bigl\{(x,T^mx,T^nx,T^{m+n}x,T^px,T^{m+p}x,T^{n+p}x,T^{m+n+p}x)\ ;\
x\in X,\ m,n,p\in\Z\bigr\}\ .
$$
Si n\'ecessaire on note $\CP(X)$ et $\CQ(X)$ au lieu de $\CP$ et $\CQ$, 
respectivement.
\end{definition}

Les relations $\RP,\RPD$ et $\RPDS$ sont li\'ees aux ensembles $\CP$ et 
$\CQ$ par le lemme facile suivant.
\begin{lemme}
\label{lem:RPRQ}
Supposons que $(X,T)$ est transitif et soient $x,y\in X$.
\begin{enumerate}
\item $(x,y)\in\RP$ si et seulement si il existe $a\in X$ avec 
$(x,y,a,a)\in\CP$.
\item
$(x,y)\in\RPD$ si et seulement si il existe $a,b,c\in X$ avec
$(x,y,a,a,b,b,c,c)\in\CQ$.
\item
$(x,y)\in\RPDS$ si et seulement si $(x,y,y,y,y,y,y,y)\in\CQ$.
\end{enumerate}
\end{lemme}
\begin{remarque}
Si le syst\`eme est minimal on montre facilement que lorsque $x,y$ sont 
r\'egionalement proximaux on a $(x,y,a,a)\in\CP$ pour tout $a\in X$.
\end{remarque}
\begin{proof}
Nous montrons seulement la deuxi\`eme affirmation, les preuves des deux 
autres \'etant similaires et plus simples. 

S'il existe $a,b,c$ tels que $(x,y,a,a,b,b,c,c)\in\CQ$, les d\'efinitions 
entra\^{\i}nent imm\'ediatement que $(x,y)\in \RPD$.

R\'eciproquement, soit $(x,y)\in\RPD$. 
Pour tout entier $i\geq 1$ soient $x'_i,y'_i, m_i$ et $n_i$ comme dans la 
d\'efinition~\ref{def:RPdeux} o\`u on a pris $\epsilon=1/i$.

Soit $u\in X$ un point d'orbite dense. 
Pour tout $i$, par densit\'e de $\{T^ku\;;\;k\in \Z\}$ et 
par continuit\'e des applications $T^{m_i}, T^{n_i}$ et $T^{m_i+n_i}$ il existe 
deux entiers $k_i$ et $\ell_i$ avec 
\begin{gather*}
d(T^{k_i}u,x)<2/i\ ;\ d(T^{\ell_i} u,y)<2/i\ ;\ 
d(T^{k_i+m_i}u,T^{\ell_i+m_i}u)<2/i\ ;\\
 d(T^{k_i+n_i}u,T^{\ell_i+n_i}u)<2/i
\text{ et }d(T^{k_i+m_i+n_i}u,T^{\ell_i+m_i+n_i}u)<2/i\ .
\end{gather*} 
En passant \`a des sous-suites nous pouvons supposer qu'il existe 
$a,b,c\in X$ tels que, quand $i\to+\infty$,
\begin{align*}
& T^{k_i+m_i}u\to a,\ T^{k_i+n_i}u\to b\text{ et } T^{k_i+m_i+n_i}u\to 
 c\\
\text{d'o\`u }\
& T^{\ell_i+m_i}u\to a,\ T^{\ell_i+n_i}u\to b\text{ et }
 T^{\ell_i+m_i+n_i}u\to  c\ .
\end{align*}
Pour tout $i$, en posant $p_i=\ell_i-k_i$ nous remarquons que
$$
(T^{k_i}u,T^{k_i+m_i}u,T^{k_i+n_i}u,T^{k_i+m_i+n_i}u,
T^{\ell_i}u,T^{\ell_i+m_i}u,T^{\ell_i+n_i}u,T^{\ell_i+m_i+n_i}u)\in\CQ
$$
et, comme $\CQ$ est ferm\'e dans $X^8$, nous avons
$(x,y,a,a,b,b,c,c)\in\CQ$.
\end{proof}

\subsection{Deux exemples}
Des exemples plus significatifs (syst\`emes distaux, nilsyst\`emes 
d'ordre $2$) seront \'etudi\'es dans les sections 
suivantes.
\subsubsection{Rotations minimales}
Soit $(X,T)$ une rotation minimale. On rappelle que $X$ est un groupe 
ab\'elien compact et que $T$ est la multiplication par une constante.
On v\'erifie directement sur les 
d\'efinitions que $\CP$ et $\CQ$ co•ncident avec les ensembles de 
parall\'elogrames et  parall\'el\'epip\`edes  d\'efinis dans les sections~2.3 
et~3.5 de~\cite{HK}:
\begin{gather*}
 \CP=\bigl\{(x,sx,tx,stx)\;;\;x,s,t\in X\bigr\}\ ;\\
\CQ=\bigl\{(x,sx,tx,stx,ux,sux,stx,stux)\;;\;x,s,t,u\in X\bigr\}
\ .
\end{gather*}

\subsubsection{Syst\`emes faiblement m\'elangeants}
On rappelle qu'un syst\`eme minimal $(X,T)$ est 
\emph{faiblement m\'elangeant} si son carr\'e cart\'esien est transitif 
et il s'ensuit que toute puissance 
cart\'esienne de $(X,T)$ est transitive. 
On en d\'eduit que, pour tout 
$k\geq 1$,  tout ouvert non vide de $X^k$ invariant par 
$T\times\dots\times T$ est dense dans $X^k$.

Montrons que $\CP=X^4$. Il suffit de montrer que qulelque soient les 
ouverts non vides  $U_0,\dots,U_3$ de $X$ il existe un parall\'elogramme 
dans $U_0\times\dots\times U_3$. Comme $X$ est faiblement m\'elangeant il
 existe un entier $n$ tel que
$(T\times T)^n(U_0\times U_1)\cap(U_2\times U_3)
\neq\emptyset$.
Il existe donc $x\in U_0$  et $y\in U_1$ avec $T^nx\in U_2$ et 
$T^ny\in U_3$. Par continuit\'e de $T^n$ et minimalit\'e, 
il existe un entier $m$ avec $T^mx\in U_1$ et $T^{m+n}x\in U_3$. 
Ainsi, $(x,T^mx,T^nx,T^{m+n}x)$ appartient \`a 
$\CP\cap(U_0\times\dots\times U_3)$ qui est donc non vide. Notre 
affirmation est d\'emontr\'ee.

Montrons maintenant que $\CQ=X^8$.  Soient $U_0,\dots,U_7$ huit 
ouverts non vides de $X$. Il existe un entier $p$ tel que 
$T^pU_i\cap U_{i+4}\neq\emptyset$ pour $0\leq i\leq 3$. D'apr\`es ce qui 
pr\'ec\`ede il existe $x\in X$ et $m,n\in\Z$ avec 
$x\in U_0\cap T^{-p}U_4$, \dots, $T^{m+n}x\in U_3\cap T^{-p}U_7$.
On obtient ainsi un parall\'el\'epip\`ede dans $U_0\times\dots\times U_7$, 
et notre affirmation est d\'emontr\'ee.

D'apr\`es le lemme~\ref{lem:RPRQ}, les relations $\RP,\RPD$ et $\RPDS$ 
sont toutes triviales (leur graphe est $X\times X$).

\subsection{Quelques questions ouvertes}

Nous ne savons pas si dans le cas g\'{e}n\'{e}ral la relation de double 
proximalit\'e r\'egionale est une 
relation d'\'equivalence, ni si le quotient de $X$ par cette relation 
est le nilfacteur d'ordre deux maximal. Ces questions seraient plus 
faciles \`a r\'esoudre si on pouvait montrer que  
cette relation passe bien aux facteurs (voir la remarque apr\`es le 
lemme~\ref{lem:facteur}).

Nous ne savons pas non plus si l'ensemble $\CP$ est li\'e de fa\c con 
simple au facteur \'equicontinu maximal (voir la remarque apre le 
lemme~\ref{lem:Ppi}).
On remarque facilement que $\CP$ et $\CQ$ v\'erifient une grande 
partie des propri\'et\'es de structures 
de parall\'elogrammes et parall\'el\'epip\`edes d\'efinies dans~\cite{HK}, au 
moins dans le cas o\`u le syst\`eme est minimal.
Cependant, nous ne savons pas si dans le cas 
 g\'en\'eral $\CP$ et $\CQ$ v\'erifient les 
propri\'et\'es de  \guillemotleft~transitivit\'e~\guillemotright, 
ni non plus si $\CQ$ v\'erifie en g\'en\'eral la 
propri\'et\'e de \guillemotleft~fermeture des 
parall\'el\'epip\`edes~\guillemotright\ introduites 
dans~\cite{HK}.

Cependant toutes ces 
propri\'et\'es sont satisfaites dans le cas des syst\`emes 
distaux(voir la d\'efinition dans la 
section~\ref{subsec:def-distal}), ce qui permet d'appliquer 
dans ce cas les r\'esultats de~\cite{HK}.

\section{Pr\'eliminaires}

\subsection{Quelques propri\'et\'es de $\CP$ et de $\CQ$}
\label{subsec:PQ}

Le lemme suivant est une application imm\'ediate des d\'efinitions.
 Les termes 
\guillemotleft~permutation euclidienne~\guillemotright\ et \guillemotleft~face~\guillemotright\ sont d\'efinis 
dans~\cite{HK}.

\begin{lemme}\label{lem:premierespropri}\strut
\begin{enumerate}
\item
$\CP$ est invariant sous les permutations euclidiennes du carr\'e.

\item
$\CQ$ est invariant sous les permutations euclidiennes du cube. 
\item
Pour tout $\bx\in\CP$, $(\bx,\bx)\in\CQ$.
\item
Si $\ubx\in\CQ$, toute face de $\ubx$ appartient \`a $\CP$.
\end{enumerate}
Si de plus $(X,T)$ est minimal nous avons :
\begin{enumerate}
\item
Pour tout $a,b\in X$, $(a,b,a,b)\in\CP$.
\item 
Soient $x_0,x_1,x_2\in X$. Il existe $x_3\in X$ avec 
$(x_0,x_1,x_2,x_3)\in\CP$.
\end{enumerate}
\end{lemme}

\begin{lemme}
\label{lem:Ppi}
Supposons que le syst\`eme $(X,T)$ est minimal.
On rappelle que $\pi\colon X\to Z$ est la projection sur le facteur 
\'equicontinu maximal.
 Pour tout $\bx=(x_0,x_1,x_2,x_3)\in\CP$ on a
$\pi(x_0)\pi(x_1)\inv\pi(x_2)\inv\pi(x_3)=1$.
\end{lemme}
\begin{remarque}
Nous ne savons pas si la r\'eciproque de ce lemme est vraie en 
g\'en\'eral.
\end{remarque}
\begin{proof}
En effet, cette propri\'et\'e est vraie lorsque 
$\bx=(x,T^mx,T^nx,T^{m+n}x)$ pour un certain $x\in X$ et certains 
$m,n\in\Z$. Elle est donc vraie pour tout parall\'elogramme par densit\'e.
\end{proof}

\begin{lemme}
\label{lem:facteur}
Soient $(Y,S)$ un facteur de $(X,T)$ et $\phi\colon X\to Y$ 
l'application facteur. Alors :
\begin{enumerate}
\item
$\CP(Y)$ est l'image de $\CP(X)$ par l'application $\phi\type 
2:=\phi\times\phi\times\phi\times\phi\colon X^4\to Y^4$;
\item
$\CQ(Y)$ est l'image de $\CQ(X)$ par l'application $\phi\type 
3:=\phi\times\dots\times\phi\colon X^8\to Y^8$;
\item
Si $x,y\in X$ sont r\'egionalement proximaux alors $\phi(x)$ et 
$\phi(y)$ sont r\'egionalement proximaux;
\item
si $x,y\in X$ sont bi-r\'egionalement proximaux (resp. fortement 
bi-r\'egionalement proximaux) alors $\phi(x)$ et 
$\phi(y)$ sont bi-r\'egionalement proximaux (resp. fortement 
bi-r\'egionalement proximaux).
\end{enumerate}
\end{lemme}
\begin{remarque}
En fait, on peut montrer que $\RP(Y)$ est l'image de $\RP(X)$ par l'application 
$\phi\times\phi$. Nous ne savons pas si la propri\'et\'e analogue est 
vraie pour la relation de bi-proximalit\'e r\'egionale.
\end{remarque}

\begin{proof} Les deux premi\`eres affirmations sont imm\'ediates d'apr\`es 
les d\'efinitions. Les deux suivantes s'en d\'eduisent imm\'ediatement gr\^ace 
au lemme~\ref{lem:RPRQ}.
\end{proof}

\subsection{Compl\'{e}ments sur la relation de proximalit\'e r\'egionale}

Le r\'{e}sultat suivant appara\^\i t dans~\cite{Aus} comme un corollaire 
de la d\'emonstration du th\'eor\`eme~\ref{th:Aus}.
\begin{corollaire*}[\cite{Aus}, Chapitre 9, corollaire 10]
Supposons que le syst\`eme $(X,T)$ est minimal.
Soit $(x,y)\in\RP$. Pour tout $\epsilon>0$ il existe $x'\in X$ et 
$n\in Z$ avec 
$$
 d(x',x)<\epsilon,\ d(T^nx',x)<\epsilon\text{ et }d(T^ny,x)<\epsilon\ .
$$
\end{corollaire*}
On en d\'eduit facilement :
\begin{corollaire}
\label{cor:RP-distaux}Supposons que le syst\`eme $(X,T)$ est minimal.
Soient $(x,y)\in\RP$. Pour tout $\epsilon>0$ et toute famille finie 
$(z_1,\dots,z_k)$ de points de $X$ il existe $z'_1,\dots,z'_k\in X$ 
et $n\in\Z$ avec
$$
d(T^ny,x)<\epsilon 
\text{ et, pour tout }i\in\{1,\dots,k\},\ d(z'_i,z_i)<\epsilon
\text{ et }
d(T^nz'_i,z_i)<\epsilon
$$
\end{corollaire}
\begin{proof}
Par minimalit\'e, il existe pour tout $i$ un entier $k_i$ avec 
$x\in T^{k_i}B(z_i,\epsilon)$. Choisissons $\delta$ avec 
$0<\delta<\epsilon$ et $B(x,\delta)\subset  T^{k_i}B(z_i,\epsilon)$ 
pour tout $i$. D'apr\`es le corollaire pr\'ec\'edent, il existe $x'\in 
B(x,\delta)$ et $n\in\Z$ avec $d(T^ny,x)<\epsilon$ et 
$d(T^nx',x)<\delta$. L'entier $n$ et les points $z'_i=T^{-k_i}x'$ v\'erifient la 
condition de l'\'enonc\'e. 
\end{proof}

\subsection{Couples proximaux, syst\`emes distaux}
\label{subsec:distal1} 
Nous rappelons la d\'efinition et les premi\`eres propri\'et\'es des syst\`emes 
distaux.

\begin{definition}
\label{def:distal1}
 Soit $(X,T)$ un syst\`eme.
On dit que le couple $(x,y)\in X\times X$ est \emph{distal} si 
$$
 \inf_{\bn\in\Z}d(T^\bn x,T^n y)>0
$$
et dans le cas contraire on dit que $(x,y)$ est \emph{proximal}.

On dit que le syst\`eme $(X,T)$ est 
\emph{distal} si tout couple $(x,y)\in X\times X$ avec $x\neq 
y$ est distal.
\end{definition}

\begin{proposition}[\cite{Aus}, chapitres 5 et 7]
\label{prop:aus-distaux}\strut
\begin{enumerate}
\item
Le produit cart\'esien d'une famille finie de 
syst\`emes distaux est un syst\`eme distal.
\item
$(X,T)$ est un syst\`eme distal et si 
$Y$ est un sous-ensemble ferm\'e de $X$ invariant par $T$, 
alors $(Y,T)$ est distal.
\item
\label{it:aus-distaux2}
Si un syst\`eme  est distal et transitif  alors il est minimal.
\item\label{it:aus-distaux1}
Tout facteur d'un syst\`eme distal est distal.
\item\label{it:aus-distaux3}
Soient $(X,T)$ un syst\`eme distal et $p\colon X\to Y$ un facteur. 
Supposons que $Y$ est minimal. Alors $p$ est une application ouverte.
\end{enumerate}
\end{proposition}

La propri\'et\'e~\ref{it:aus-distaux3}
 est d\'emontr\'ee dans~\cite{Aus} (chapitre 7,  th\'eor\`eme 3) sous l'hypoth\`ese 
suppl\'{e}mentaire que $X$ est minimal, mais il est tr\`es facile d'en 
d\'eduire le cas g\'en\'eral.
\begin{proof}
Soient $x\in X$ et $V$ un voisinage de $x$ dans $X$. Nous devons 
montrer que $p(V)$ est un voisinage de $p(x)$ dans $Y$.

Soit $X'$ l'orbite ferm\'ee de $x$ dans $X$. Muni de la restriction de
$T$,  $X'$ est un syst\`eme distal et transitif, donc 
minimal d'apr\`es~\ref{it:aus-distaux2}.
La restriction de $p$ \`a $X'$ est un facteur. Ainsi, $p(V\cap X')$ est 
un voisinage de $p(x)$ dans $Y$, donc \'egalement $p(V)$.
\end{proof}

\begin{proposition}
\label{prop:distalRP}
Soit $(X,T)$ un syst\`eme transitif tel que la relation $\RPDS$ soit 
l'\'egalit\'e. Alors ce syst\`eme est distal et minimal.
\end{proposition}
La conclusion reste clairement valable si on suppose que l'une 
des relations $\RP$ ou $\RPD$ est l'\'egalit\'e, puisque ces deux 
relations sont moins fines que la relation $\RPDS$.
Pour montrer le th\'eor\`eme~\ref{th:main} on peut donc sans perdre en 
g\'en\'eralit\'e se restreindre au cas des syst\`emes distaux minimaux.

Nous commen\c cons par un lemme.
\begin{lemme}
\label{lem:distalRP}
Soient $(X,T)$ un syst\`eme et $(x,y)$ un couple proximal. Supposons 
que l'orbite ferm\'ee de $y$ est minimale sous la transformation $T$. 
Alors $(x,y)\in\RPDS$.
\end{lemme}
\begin{remarque}
Un point dont l'orbite ferm\'ee est minimale est appel\'e \guillemotleft~almost 
periodic~\guillemotright\ 
dans~\cite{Aus}, mais nous pr\'ef\'erons \'eviter ce terme \`a cause du 
risque de confusion avec la notion de \guillemotleft~suite presque p\'eriodique~\guillemotright\ 
consid\'er\'ee dans la section~\ref{sec:nilsuites}.
\end{remarque}
\begin{proof}[D\'emonstration du lemme~\ref{lem:distalRP}]
Par hypoth\`ese il existe une suite d'entiers $(k_i)$ telle que 
$d(T^{k_i}x,T^{k_i}y)\to 0$. Quitte \`a remplacer cette suite par une 
sous-suite on peut supposer que les  suites $(T^{k_i}x)$ et 
$(T^{k_i}y)$ convergent toutes deux vers le m\^eme point $z\in X$. 
Comme $z$ appartient \`a l'orbite ferm\'ee de $y$ et 
que cette orbite est minimale, $y$ appartient \`a l'orbite ferm\'ee de $z$

Soit $\epsilon>0$.  Il
 existe un entier $\ell$ avec $d(T^\ell z,y)<\epsilon/3$. 
Par continuit\'e de $T^\ell$, en posant $m=k_i+\ell$ pour $i$ assez 
grand, on a 
$$
 d(T^mx,y)<\epsilon/2\text{ et }d(T^my,y)<\epsilon/2\ .
$$
Choisissons $\delta$ avec $0<\delta<\epsilon$ tel que 
$d(T^mw,T^my)<\epsilon/2$ pour tout $w\in B(y,\delta)$. En proc\'edant 
comme pr\'{e}c\'{e}demment on obtient un entier $n$ tel que
\begin{align*}
&  d(T^nx,y)<\delta\text{ et }d(T^ny,y)<\delta\\
\text{d'o\`u } \qquad
& d(T^{m+n}x,T^my)<\epsilon/2\text{ et } d(T^{m+n}y,T^my)<\epsilon/2
\\
\text{et enfin } \qquad
& d(T^{m+n}x,y)<\epsilon\text{ et } d(T^{m+n}y,y)<\epsilon
\end{align*}
ce qui ach\`eve la d\'emonstration.
\end{proof}

\begin{proof}[D\'emonstration de la proposition~\ref{prop:distalRP}]
Comme $(X,T)$ est transitif, il existe un point $x\in X$ dont l'orbite est 
dense. D'apr\`es le th\'eor\`eme~10 du chapitre~6 de~\cite{Aus} il existe un 
point $y\in X$ dont l'orbite ferm\'ee soit minimale  et  tel que 
$(x,y)$ soit un couple proximal. D'apr\`es le lemme~\ref{lem:distalRP}, 
$(x,y)\in\RPDS$ donc $x=y$ par hypoth\`ese. L'orbite ferm\'ee de $y$ est 
donc \'egale \`a $X$ et $(X,T)$ est minimal.

Soit maintenant $(w,z)$ un couple proximal. Comme $X$ est minimal, 
l'orbite ferm\'ee de $z$ est \'egale \`a $X$ et est minimale sous $T$. Le 
lemme~\ref{lem:distalRP} entra\^\i ne que $(w,z)\in\RPDS$ et on a donc 
$w=z$ par hypoth\`ese. Le syst\`eme $(X,T)$ est donc distal.
\end{proof}

\subsection{Syst\`emes distaux avec plusieurs transformations}
\label{subsec:def-distal}
Nous aurons besoin de nous placer dans un cadre plus g\'en\'eral que celui
 consid\'er\'e jusqu'ici. 

 Soient $X$ un espace compact muni d'une distance $d$, 
$k\geq 1$ un entier  et  $T_1,\dots,T_k$ des hom\'eomorphismes 
de $X$ qui commutent entre eux.
Ces transformations induisent une action 
de $\Z^k$ sur $X$ par hom\'eomorphismes :
$$
\text{pour }\bn=(n_1,\dots,n_k)\in\Z^k
\text{ on note }T^\bn=T_1^{n_1}\dots T_k^{n_k}\ .
$$
Le compact $X$ muni des transformations $T_1,\dots,T_k$ sera appel\'e un 
\emph{$\Z^k$-syst\`eme}, ou simplement un syst\`eme.

Le syst\`eme $(X,T_1,\dots,T_k)$ est dit transitif s'il existe au moins 
un point $x\in X$ dont l'orbite $\{T^\bn\;;\; \bn\in\Z^k\}$ est dense 
dans $X$, et il est dit minimal si l'orbite de tout point est dense.

Un \emph{facteur} de $(X,T_1,\dots,T_k)$ est un la donn\'ee d'un syst\`eme 
$(Y,S_1,\dots,S_k)$ et d'une application $p\colon X\to Y$, continue 
et surjective, telle que $p\circ T_i=S_i\circ p$ pour $1\leq i\leq k$.

Les d\'efinitions  de la section~\ref{subsec:distal1} s'\'etendent sans modification autre 
que de notations :

 Soit $(x,y)\in X\times X$. 
On dit que le couple $(x,y)$ est \emph{distal} si 
$$
 \inf_{\bn\in\Z^k}d(T^\bn x,T^\bn y)>0
$$
et dans le cas contraire on dit que $(x,y)$ est \emph{proximal}.
On dit que le syst\`eme $(X,T_1,\dots,T_k)$ est 
\emph{distal} si tout couple $(x,y)\in X\times X$ avec $x\neq 
y$ est distal.

\begin{proposition}[\cite{Aus}, chapitres 5 et 7]
\label{prop:aus-distaux2}
Les r\'esultats de la proposition~\ref{prop:aus-distaux}  sont valables 
pour les $\Z^k$-syst\`emes.
\end{proposition}

\subsection{Premi\`eres propri\'et\'es de $\CP$ et $\CQ$ dans un syst\`eme 
minimal distal}

\begin{notations}
Soit $(X,T)$ un syst\`eme.
On note $T\type 2, T\type 2_1$ et $T\type 2_2$ les transformations de 
$X^4$ donn\'ees par :\\ 
\null\quad pour $\bx=(x_0,x_1,x_2,x_3)$,
\begin{align*}
T\type 2\bx& =(Tx_0,Tx_1,Tx_2,Tx_3)\ ;\\
T\type 2_1\bx &= (x_0,x_1,Tx_2,Tx_3)\ ;\\
T\type 2_2\bx &=(x_0,Tx_1,x_2,Tx_3)\ .
\end{align*}
On note $T\type 3,T\type 3_1,T\type 3_2$ et $T\type 3_3$ les 
transformations de $X\type 8$ donn\'ees par :\\
\null\quad pour $\ubx=(x_0,x_1,\dots,x_7)$,
\begin{align*}
T\type 3\ubx &= (Tx_0,Tx_1,Tx_2,Tx_3,Tx_4,Tx_5,Tx_6,Tx_7)\ ;\\
T\type 3_1\ubx &= (x_0,x_1,x_2,x_3,Tx_4,Tx_5,Tx_6,Tx_7)\ ;\\
T\type 3_2\ubx &= (x_0,x_1,Tx_2,Tx_3,x_4,x_5,Tx_6,Tx_7)\ ;\\
T\type 3_3\ubx &= (x_0,Tx_1,x_2,Tx_3,x_4,Tx_5,x_6,Tx_7)\ .
\end{align*}
\end{notations}
Par construction, $\CP$ est invariant par les transformations $T\type 
2, T\type 2_1$ et $T\type 2_2$, et $\CQ$ est invariant par les 
transformations $T\type 3,T\type 3_1,T\type 3_2$ et $T\type 3_3$.
\begin{lemme}
\label{lem:Qtrnasitif}
Soit $(X,T)$ un syst\`eme distal minimal. 
Alors $(\CP,T\type 2,T\type 2_1,T\type 2_2)$ et 
$(\CQ,T\type 3,T\type 3_1,T\type 3_2,T\type 3_3)$ sont des syst\`emes 
distaux et minimaux.
\end{lemme}

\begin{proof}
Nous montrons seulement la deuxi\`eme affirmation. 

Comme $(X,T)$ est distal, il est imm\'ediat que 
$X^8$ muni des transformations $T\type 3,T\type 3_1,T\type 3_2$ et $T\type 3_3$ 
est distal. 
Comme $\CQ$ est un ferm\'e de $X^8$ invariant par ces transformations, 
$(\CQ,T\type 3,T\type 3_1,T\type 3_2,T\type 3_3)$ est distal. Pour 
montrer que ce syst\`eme est minimal il suffit donc de montrer qu'il 
est transitif.

Soit $z\in X$ et posons $\ubz=(z,z,\dots,z)\in X^8$. Par 
d\'efinition, $\ubz\in\CQ$. Nous montrons que l'orbite de $\ubz$ est 
dense dans $\CQ$.

 Soit $\ubx=(x_0,\dots,x_7)\in\CQ$ et $\epsilon>0$. Par 
construction, il existe $x\in X$ et des entiers $m,n,p$ tels que
chaque sommet de $\ubx$ soit \`a une distance $<\epsilon/2$ du sommet 
correspondant de $(x,T^mx,T^nx,\dots,T^{m+n+p}x)$.
Par minimalit\'e et continuit\'e des applications $T^m, T^n$ et $T^{m+n}$,
 il existe un entier $k$ tel que $T^kz$ soit 
suffisamment voisin de $x$ pour que
$$
 d(T^kz,x_0)<\epsilon,\ 
d(T^{k+m}z,x_1)<\epsilon,\ d(T^{k+n}z,x_2)<\epsilon,\ \dots,\ 
d(T^{k+m+n+p}z,x_7)<\epsilon\ .
$$
Ainsi, chaque sommet de ${T\type 3}^k{T\type 3_1}^m{T\type 
k_2}^n{T\type 3_3}^n\ubz$ est \`a une distance $<\epsilon$ du 
sommet correspondant de $\ubx$ 
\end{proof}

\section{Le cas des syst\`emes distaux}
\label{sec:distaux}

Nous continuons ici l'\'etude des parall\'elogrammes et 
parall\'{e}l\'{e}pip\`{e}des dans les syst\`emes distaux, le r\'esultat principal 
\'etant le th\'eor\`eme~\ref{th:struc-distal} qui permet d'utiliser dans la 
suite la machinerie d\'evelopp\'ee dans~\cite{HK}.

Dans toute  cette section, $(X,T)$ est un syst\`eme minimal 
distal. On rappelle que $\pi\colon X\to Z$ est le facteur \'equicontinu 
maximal de $X$.

\subsection{Parall\'elogrammes dans un syst\`eme distal}

\begin{theoreme}
\label{th:para-distal}
Soit $\bx=(x_0,x_1,x_2,x_3)\in X^4$. 
\begin{enumerate}
\item \label{it:para-distal1}
Pour que $\bx$ appartienne \`a $\CP$ il faut et il suffit que 
$$
\pi(x_0)\pi(x_1)\inv\pi(x_2)\inv\pi(x_3)=1\ .
$$
\item\label{it:para-distal2}
Pour que $\bx$ appartienne \`a $\CP$ il faut et il suffit qu'il existe deux suites d'entiers $(m_i)$ 
et $(n_i)$ avec
$$
 T^{m_i}x_0\to x_1,\;T^{n_i}x_0\to x_2\text{ et }T^{m_i+n_i}x_0\to 
x_3\ .
$$
\end{enumerate}
\end{theoreme}
Ainsi, $\CP$ est la \guillemotleft~structure de parall\'elogrammes~\guillemotright\ sur $X$ associ\'ee \`a 
la projection $\pi\colon X\to Z$ (\cite{HK}).

Nous commen\c cons par un lemme.
\begin{lemme}
\label{lem:para-distal1}
Soient $x_0\in X$ et  $K$ l'orbite ferm\'ee  dans $X^4$ de 
$(x_0,x_0,x_0,x_0)$ pour les transformations $T\type 2_1$ et $T\type 
2_2$. 

Si $(x_0,x'_1,x_2,x_3)\in K$ et si $(x'_1,x_1)\in\RP$ 
alors $(x_0,x_1,x_2,x_3)\in K$.
\end{lemme}
\begin{proof}[D\'emonstration du lemme~\ref{lem:para-distal1}]
Soit $\epsilon >0$. Notons 
$$
 Y=\{(y_2,y_3)\in X^2\;;\; (x_0,x'_1,y_2,y_3)\in K\}\ .
$$
Alors $K$ est un ferm\'e de $X\times X$ contenant $(x_2,x_3)$. 
Comme $K$ est invariant par 
$T\type 2_2$, $Y$ est invariant par $T\times T$. 
Munissons $Y$ de cette transformation. Alors l'application 
$(y_2,y_3)\mapsto y_3$ est un facteur de $(Y,T\times T)$ dans 
$(X,T)$. D'apr\`es la proposition~\ref{prop:aus-distaux2}, cette application 
est ouverte. Il existe donc $\delta>0$ tel que 
$$
 \text{pour tout }y_3\in B(x_3,\delta)\text{ il existe } y_2\in 
B(x_2,\epsilon)\text{ avec } (x_0,x'_1,y_2,y_3)\in K\ .
$$

Comme $x'_1$ et $x_1$ sont r\'egionalement proximaux, d'apr\`es
le corollaire~\ref{cor:RP-distaux}  il existe $y_3\in X$ 
et $n\in\Z$ avec
$$
 d(y_3,x_3)<\delta,\ d(T^ny_3, x_3)<\epsilon\text{ et 
}d(T^nx'_1,x_1)<\epsilon\ .
$$
D'apr\`es le choix de $\delta$, il existe $y_2\in X$ avec 
$d(y_2,x_2)<\epsilon$ et $(x_0,x'_1,y_2,y_3)\in K$. 
Comme $K$ est invariant par $T\type 2_1$, 
$(x_0,T^nx'_1,y_2,T^ny_3)\in K$. 

En faisant tendre $\epsilon$ vers $0$ on obtient que 
$(x_0,x_1,x_2,x_3)\in K$.
\end{proof}

\begin{proof}[D\'emonstration du th\'eor\`eme~\ref{th:para-distal}]
La propri\'et\'e~\ref{it:para-distal2} signifie que $\bx$ appartient \`a 
l'ensemble $K$ introduit dans le lemme~\ref{lem:para-distal1}. 
Cette propri\'et\'e entra\^\i ne imm\'ediatement que $\bx$ appartient \`a $\CP$, 
ce qui d'apr\`es le lemme~\ref{lem:Ppi} entra\^\i ne la propri\'et\'e~\ref{it:para-distal1}.

 Supposons que $\bx$ v\'erifie la propri\'et\'e de~\ref{it:para-distal1}.

 Choisissons deux suites 
d'entiers $(n_i)$  et $(k_i)$ telles que $T^{n_i}x_0\to x_2$ et 
$T^{k_i}x_0\to x_3$ et posons $m_i=k_i-n_i$. En rempla\c cant les suites 
par des sous-suites on peut supposer que la suite $(T^{m_i}x_0)$ 
converge dans $X$ vers un point $x'_1$.

Alors par construction $(x_0,x'_1,x_2,x_3)$ appartient \`a l'ensemble 
$K$ d\'efini dans le lemme~\ref{lem:para-distal1}. En particulier,
 $(x_0,x'_1,x_2,x_3)$ appartient \`a $\CP$ et donc
$\pi(x_0)\pi(x'_1)\inv$ $\pi(x_2)\inv\pi(x_3)=1$. L'hypoth\`ese entra\^{\i}ne 
alors que 
$\pi(x'_1)=\pi(x_1)$ et donc que $(x'_1,x_1)\in\RP$. 
D'apr\`es le lemme~\ref{lem:para-distal1}, $(x_0,x_1,x_2,x_3)$ 
appartient \`a $K$. La propri\'et\'e~\ref{it:para-distal2} est satisfaite.
\end{proof}

\subsection{Parall\'el\'epip\`edes dans un syst\`eme distal}
Dans cette section nous supposons encore que $(X,T)$ est un syst\`eme 
minimal distal et
nous montrons :
\begin{theoreme}\label{th:struc-distal}
$(\CP,\CQ)$ est une structure de parall\'el\'epip\`edes sur $X$ au sens 
de~\cite{HK}.
\end{theoreme}

Nous identifions $X^8$ et $X^4\times X^4$ : pour 
$\bx=(x_0,x_1,x_2,x_3)$ et $\by=(y_0,y_1,y_2,y_3)\in X^4$,
$(\bx,\by)=(x_0,\dots,x_3,y_0,\dots,y_3)\in X^8$. Commen\c cons par 
montrer :
\begin{proposition}
\label{prop:trans-Q}
$\CQ$ v\'erifie la \guillemotleft~propri\'et\'e de transitivit\'e~\guillemotright:

Soient $\bu,\bv$ et $\bw\in\CP$. Si $(\bu,\bv)$ et $(\bv,\bw)$ 
appartiennent \`a $\CQ$ alors $(\bu,\bw)$ appartient \`a $\CQ$.
\end{proposition}

\begin{proof} [D\'emonstration de la proposition~\ref{prop:trans-Q}]
Soit $z\in X$. Notons $\bz=(z,z,z,z)\in X^4$ et 
$\ubz=(\bz,\bz)=(z,z,\dots,z)\in X^8$. On a $\ubz\in\CQ$.
D'apr\`es le lemme~\ref{lem:Qtrnasitif}, $(\CQ,T\type 3,T\type 3_1,T\type 
3_2,T\type 3_3)$ est minimal et donc $\CQ$ est l'orbite ferm\'ee de 
$\ubz$ pour ces quatre transformations.

Il existe donc quatre suites d'entiers $(m_i)$, $(n_i)$, $(p_i)$ et 
$(q_i)$ telles que
$$
 {T\type 3_1}^{m_i} {T\type 3_2}^{n_i} {T\type 3_3}^{p_i}
{T\type 3}^{q_i}\,(\bu,\bv)\to\ubz\text{ dans }X^8\ .
$$
Ainsi,
$$
 {T\type 2_1}^{m_i}{T\type 2_2}^{n_i}{T\type 2}^{q_i}\bu\to\bz
\text{ et }
 {T\type 2_1}^{m_i}{T\type 2_2}^{n_i}{T\type 2}^{p_i+q_i}\bv\to\bz
\text{ dans }X^4\ .
$$
Quitte \`a remplacer les suites par des sous-suites, on peut supposer 
que la suite
$$
  {T\type 2_1}^{m_i}{T\type 2_2}^{n_i}{T\type 2}^{p_i+q_i}\bw
$$
converge  dans $X^4$ vers un point $\bj\in\CP$.

Pour tout $i$,
$$
 \bigl( {T\type 2_1}^{m_i}{T\type 2_2}^{n_i}{T\type 2}^{p_i+q_i}\bv,
{T\type 2_1}^{m_i}{T\type 2_2}^{n_i}{T\type 2}^{p_i+q_i}\bw\bigr)
={T\type 3_1}^{m_i}{T\type 3_2}^{n_i}{T\type 3}^{p_i+q_i}\,(\bv,\bw)
$$
qui appartient \`a $\CQ$, donc \`a la limite $(\bz,\bj)\in\CQ$.
Par ailleurs, pour tout $i$,
$$
 {T\type 3_1}^{m_i}{T\type 3_2}^{n_i}{T\type 3_3}^{p_i}{T\type 3}^{q_i}
\,(\bu,\bw)=
\bigl( {T\type 2_1}^{m_i}{T\type 2_2}^{n_i}{T\type 2}^{q_i}\,\bu,
{T\type 2_1}^{m_i}{T\type 2_2}^{n_i}{T\type 2}^{p_i+q_i}\bw\bigr)
$$
qui converge vers $(\bz,\bj)$ dans $X^8$. Ainsi, $(\bz,\bj)$ 
appartient \`a l'orbite ferm\'ee de $(\bu,\bw)$ sous les transformations
$T\type 3_1,T\type 3_2,T\type 3_3$ et $T\type 3$. 

Comme cette orbite 
est minimale, $(\bu,\bw)$ appartient \`a l'orbite ferm\'ee de $(\bz,\bj)$ 
sous ces transformations. Comme $(\bz,\bj)$ appartient \`a $\CQ$, 
$(\bu,\bw)$ appartient \`a $\CQ$.
\end{proof}

Le th\'eor\`eme~\ref{th:struc-distal} se d\'eduit maintenant imm\'ediatement 
de la proposition suivante.
\begin{proposition}
\label{prop:EFFP}
Pour $(x_0,\dots,x_6)\in X^7$  les propri\'et\'es suivantes sont 
\'equivalentes.
\begin{enumerate}
\item\label{it:EFPP1}
 $(x_0,x_1,x_2,x_3)$, $(x_0,x_1,x_4,x_5)$ et
$(x_0,x_2,x_4,x_6)$ appartiennent \`a $\CP$;
\item\label{it:EFPP2}
il existe $x_7\in X$ avec $(x_0,\dots,x_6,x_7)\in\CQ$;
\item\label{it:EFPP3}
Il existe trois suites d'entiers 
 $(m_i)$, $(n_i)$ et $(p_i)$ avec
$ T^{m_i}x_0\to x_1$, $T^{n_i}x_0\to x_2$, $T^{m_i+n_i}x_0\to x_3$, 
 $T^{p_i}x_0\to x_4$, $T^{m_i+p_i}x_0\to x_5$ et  
$T^{n_i+p_i}x_0\to x_6$.
\end{enumerate}
\end{proposition}

Pour montrer cette proposition nous aurons besoin de deux lemmes. 

\begin{lemme}
\label{lem:Q6}
Soient $x_0\in X$ et $L$ l'adh\'erence dans $X^5$ de l'ensemble
$$
 \bigl\{(T^mx_0,T^nx_0,T^{m+n}x_0,T^px_0,T^{p+m}x_0)\;;\;
m,n,p\in\Z\bigr\}\ .
$$
Soit $(x_1,\dots,x_5)\in K$ et $x'_2,x'_4\in X$ tels que
$(x_2,x'_2)\in\RP$ et $(x_4,x'_4)\in\RP$.  
Alors  $(x_1,x'_2,x_3,x'_4,x_5)\in L$.
\end{lemme}
\begin{proof}[D\'emonstration du lemme~\ref{lem:Q6}]
Montrons que  $(x_1,x'_2,x_3,x_4,x_5)\in L$. Soit $\epsilon>0$.

Posons 
$$
 F=\bigl\{ (y_1,y_3,y_5)\in X^3\;;\; (y_1,x_2,y_3,x_4,y_5)\in 
L\bigr\}\ .
$$
Alors $F$ est invariant par la transformation $T\times T\times T$ et 
$(F,T\times T\times T)$ est distal. L'application 
$(y_1,y_3,y_5)\mapsto y_3$ est un facteur de ce syst\`eme sur $(X,T)$. 
Cette application est donc ouverte et il existe $\delta>0$ tel que, 
pour tout $y_3\in X$, si $d(x_3,y_3)<\delta$, il existe 
$y_1,y_5\in X$ avec
\begin{equation}
\label{eq:y1y5}
d(x_1,y_1)<\epsilon,\ d(x_5,y_5)<\epsilon\text{ et 
}(y_1,x_2,y_3,x_4,y_5)\in L\ .
\end{equation}
Comme $x_3$ et $x'_3$ sont r\'egionalement proximaux, il existe 
$y_3\in X$ et $n\in\Z$ avec
$$
 d(x_3,y_3)<\delta,\ d(T^ny_3,x_3)<\epsilon\text{ et 
}d(T^nx_2,x'_2)<\epsilon\ .
$$
Soient $y_1,y_5$ v\'erifiant~\eqref{eq:y1y5}. Alors 
$(y_1,T^nx_2,T^ny_3, x_4,y_5)\in L$. 
 
En passant \`a la limite quand $\epsilon\to 0$
 on 
obtient que $(x_1,x'_2,x_3,x_4,x_5)\in L$.

Par la m\^eme d\'emonstration on obtient que 
$(x_1,x'_2,x_3,x'_4,x_5)$ 
appartient \`a $L$.
\end{proof}

\begin{lemme}
\label{lem:Q7}
Soient $x_0\in X$ et $K$ l'adh\'erence dans $X^6$ de l'ensemble
$$
 \bigl\{(T^mx_0,T^nx_0,T^{m+n}x_0,T^px_0,T^{p+m}x_0,T^{p+n}x_0)\;;\;
m,n,p\in\Z\bigr\}\ .
$$
Soit $(x_1,\dots,x_6)\in K$ et $x'_6\in X$ tel que $x_6$ et $x'_6$ 
soient r\'egionalement proximaux. Alors  $(x_1,\dots,x_5,x'_6)\in K$.
\end{lemme}
\begin{proof}[D\'emonstration du lemme~\ref{lem:Q7}]
Fixons $\epsilon>0$.

Munissons $X^3$ des trois transformations
$$
 S_1=\id\times T\times \id,\ S_2=\id\times \id\times T\text{ et }
S_3=T\times T\times T\ .
$$
Alors $(X^3,S_1,S_2,S_3)$ est distal et 
la minimalit\'e de $(X,T)$ entra\^\i ne que ce syst\`eme est minimal.

D'autre part, $K$ est invariant par les trois transformations :
\begin{align*}
T_1 & \colon(y_1,y_2,y_3,y_4,y_5,y_6)\mapsto  
 (Ty_1,y_2,Ty_3,y_4,Ty_5,y_6)\ ;\\
T_2 & \colon(y_1,y_2,y_3,y_4,y_5,y_6)\mapsto  
 (y_1,Ty_2,Ty_3,y_4,y_5,Ty_6)\ ;\\
T_3 & \colon(y_1,y_2,y_3,y_4,y_5,y_6) \mapsto  
(y_1,y_2,y_3,Ty_4,Ty_5,Ty_6)\ .
\end{align*}
Le syst\`eme $(K,T_1,T_2,T_3)$ est distal et transitif donc il est 
minimal.

L'application 
$p\colon(y_1,y_2,y_3,y_4,y_5,y_6)\mapsto(y_4,y_5,y_6)$ 
de $K$ dans $X^3$ est continue et commute avec les transformations, 
elle est donc surjective par minimalit\'e de $X^3$, et $p$ 
est ainsi une application facteur de $(K,T_1,T_2,T_3)$ sur 
$(X^3,S_1,S_2,S_3)$. Cette 
application est donc ouverte. Il existe donc $\delta$ tel que pour 
tous
$y_4,y_5\in X$ avec
$$
 d(x_4,y_4)<\delta\text{ et }(x_5,y_5)<\delta
$$
il existe $y_1,y_2,y_3\in X$ avec
\begin{equation}
\label{eq:y1y2y3}
d(x_1,y_1)<\epsilon,\ d(x_2,y_2)<\epsilon,\ d(x_3,y_3)<\epsilon
\text{ et }
(y_1,y_2,y_3,y_4,y_5,x_6)\in K\ .
\end{equation}
Comme $x'_6$ et $x_6$ sont r\'egionalement proximaux, d'apr\`es 
le corollaire~\ref{cor:RP-distaux} il existe $y_4,y_5\in X$ 
$n\in\Z$ avec
$$
d(x_4,y_4)<\delta,\  d(T^ny_4,x_4)<\epsilon,\
 d(x_5,y_5)<\delta,\ d(T^ny_5,x_5)<\epsilon\text{ et }
d(T^nx_6,x'_6)<\epsilon\ .
$$
Soient $y_1,y_2,y_3$ comme dans~\eqref{eq:y1y2y3}. Alors
$(y_1,y_2,y_3,T^ny_4,T^ny_5,T^nx_6)\in K$.

Par passage \`a la limite,  $(x_1,\dots,x_5,x'_6)\in K$.
\end{proof}

\begin{proof} [D\'emonstration de la proposition~\ref{prop:EFFP}]
Supposons que la propri\'et\'e~\ref{it:EFPP1} est v\'erifi\'ee.

Soient $(m_i)$, $(k_i)$ et $(\ell_i)$ trois suites d'entiers telles 
que 
 $T^{m_i}x_0\to x_1$, $T^{k_i}x_0\to x_3$ et $T^{\ell_i}\to x_5$. 
Posons $n_i=k_i-m_i$ et $p_i=\ell_i-m_i$. En rempla\c cant les suites 
donn\'ees par des sous-suites on peut supposer que les suites 
$(T^{n_i}x_0)$ et $(T^{p_i}x_0)$ convergent dans $X$; soient $x'_2$ et 
$x'_4$ leurs limites respectives.

Par construction, $(x_1,x'_2,x_3,x'_4,x_5)$ appartient \`a l'ensemble 
$L$ d\'efini dans le lemme~\ref{lem:Q6}. En particulier, 
$(x_0,x_1,x'_2,x_3)\in\CP$, et comme $(x_0,x_1,x_2,x_3)$ 
appartient aussi \`a $\CP$, $x'_2$ et $x_2$ sont r\'egionalement 
proximaux. De m\^eme, $x'_4$ et $x_4$ sont r\'egionalement proximaux. 
D'apr\`es le lemme~\ref{lem:Q6}, 
$(x_1,x_2,x_3,x_4,x_5)\in L$.

Soit $K$ l'ensemble introduit dans le lemme~\ref{lem:Q7}. Alors $L$ 
est clairement l'image de $K$ par la projection
$(y_1,y_2,y_3,y_4,y_5,y_6)\mapsto (y_1,y_2,y_3,y_4,y_5)$. Il existe 
donc $x'_6\in X$ avec $(x_1,x_2,x_3,x_4,x_5,x'_6)\in K$. 
Comme $(x_0,_2,x_4,x_6)$ et $(x_0,_2,x_4,x'_6)$ appartiennent \`a
$\CP$, $x_6$ et $x'_6$ sont r\'egionalement proximaux. D'apr\`es le 
lemme~\ref{lem:Q7},  $(x_1,x_2,x_3,x_4,x_5,x_6)$ appartient \`a $K$.
Par d\'efinition de $K$, la propri\'et\'e~\ref{it:EFPP3} est v\'erifi\'ee.

Supposons la propri\'et\'e~\ref{it:EFPP3} v\'erifi\'ee et soient
 trois suites d'entiers $(m_i)$, $(n_i)$ et $(p_i)$ comme dans cette 
propri\'et\'e. En rempla\c cant ces suites par des sous-suites 
on peut supposer que la suite $(T^{m_i+n_i+n_i}x_0)$ converge vers un 
point $x_7$. On a $(x_0,\dots,x_7)\in\CQ$, et la 
propri\'et\'e~\ref{it:EFPP2} est v\'erifi\'ee. 

Il est \'evident que la propri\'et\'e~\ref{it:EFPP2} 
entra\^\i ne la 
propri\'et\'e~\ref{it:EFPP1}.
\end{proof}
\subsection{D\'emonstration de la proposition~\ref{prop:RPDRPDS} et 
d'une partie du th\'eor\`eme~\ref{th:main-distal}}
\label{subsec:RPDRPDS}

Soit $(X,T)$ un syst\`eme distal minimal. D'apr\`es le 
th\'eor\`eme~\ref{th:struc-distal}, $\CP$ et $\CQ$ forment une structure 
de parall\'el\'epip\`edes sur $X$.  D'apr\`es le lemme~\ref{lem:RPRQ}, les 
deux relations $\RPD$ et $\RPDS$ co•ncident avec la relation
$\equivq$ introduite dans la section~3.3 de~\cite{HK}. En particulier, 
$\RPD$ est une relation d'\'equivalence, et cette relation est 
clairement ferm\'ee. Dans la section~\ref{subsec:proof-main-distal} 
on montrera  que le quotient de $X$ par cette relation est le 
nilfacteur d'ordre $2$ maximal, ce qui ach\`evera la d\'emonstration du 
th\'eor\`eme~\ref{th:main-distal}.

Nous montrons maintenant la proposition~\ref{prop:RPDRPDS}.
Si la relation $\RPD$ est l'\'egalit\'e, la relation $\RPDS$ qui est plus 
fine est \'egalement l'\'egalit\'e.
Soit maintenant $(X,T)$ un syst\`eme transitif et que la relation 
$\RPDS$ soit l'\'egalit\'e. D'apr\`es la proposition~\ref{prop:distalRP}, ce 
syst\`eme est distal et minimal et d'apr\`es ce qui pr\'ec\`ede les 
relations $\RPD$ et $\RPDS$ co•ncident. \qed

\section{Le cas des nilsyst\`emes}

\subsection{Rappels}
\label{subsec:rappels-nil}
Soient $k\geq 1$ un entier et $G$ un groupe de Lie nilpotent d'ordre 
$k$.  Soient $\Gamma$,
$X=G/\Gamma$, $t\in G$ et $T$ comme dans la 
d\'efinition~\ref{def:nilsysteme}.
Nous rappelons quelques propri\'et\'es classiques utiles dans la suite, 
dont la d\'emonstration se trouve dans~\cite{AGH}, \cite{Pa} et~\cite{Le}.

\begin{proposition}
La famille des nilsyst\`emes d'ordre $k$ est stable par passage aux 
facteurs et par produit cart\'esien.
\end{proposition}

\begin{theoreme}
\label{th:rappelsnil1}
Soit $(X,T)$ un nilsyst\`eme d'ordre $k$.
\begin{itemize}
\item
$(X,T)$ est distal. En particulier, si $(X,T)$ est transitif alors il 
est minimal.
\item
Soit $Y$ l'orbite ferm\'ee d'un point de $X$. Alors $(Y,T)$ est (isomorphe \`a) 
un nilsyst\`eme d'ordre $k$.
\end{itemize}
\end{theoreme}

Une premi\`ere r\'eduction est utile avant d'\'enoncer la propri\'et\'e suivante.

Soient $G,\Gamma$ et $t$ comme plus haut et supposons que le nilsyst\`eme 
$(X,T)$ est minimal.
Soient $G_0$ la composante connexe de l'\'el\'ement unit\'e de $G$ et $G_1$ 
le sous-groupe de $G$ engendr\'e par $G_0$ et $t$. Alors $G_1$ est un 
sous-groupe ouvert de $G$ et sa projection sur $X=G/\Gamma$ est un 
ouvert de $X$, invariant par $T$ donc \'egal \`a $X$ par minimalit\'e. 
Posons $\Gamma_1= \Gamma\cap G_1$. Alors $\Gamma_1$ est un sous-groupe 
discret et cocompact de $G_1$ et  $X$ s'identifie \`a  $G_1/\Gamma_1$. 
En rempla\c cant $G$ par $G_1$ et $\Gamma$ par $\Gamma_1$ nous pouvons 
ainsi nous restreindre sans perdre en g\'en\'eralit\'e, au cas o\`u $G$  et 
$t$
v\'erifie la propri\'et\'e :
\begin{equation}
\label{eq:H1}
\tag{H1}
\text{\em $G$ est engendr\'e par la composante connexe de son \'el\'ement 
unit\'e et $t$.}
\end{equation}

Notons $G_2$ le groupe des commutateurs de $G$. On rappelle que $G_2$ 
et $G_2\Gamma$ sont des sous-groupes ferm\'es de $G$.
L'hypoth\`ese~(H1) entra\^\i ne imm\'{e}diatement que $G_2$ est connexe. 

\begin{theoreme}[voir \cite{Le}]
\label{th:rappelsnil2}
Supposons que le nilsyst\`eme $(X,T)$ est minimal et que 
l'hypoth\`ese~\emph{(H1)} est v\'erifi\'ee. 
 Alors le facteur \'equicontinu 
maximal de  ce syst\`eme est le quotient $Z=G/G_2\Gamma$ de $X=G/\Gamma$ 
par l'action de $G_2$.
\end{theoreme}

Une deuxi\`eme r\'eduction sera utile dans la suite. Notons $\CZ(G)$ le 
centre de $G$. Alors
$Z(G)\cap\Gamma$ est un sous-groupe distingu\'e et ferm\'e
de $G$. En rempla\c cant $G$ et $\Gamma$ par leurs quotients par 
ce sous-groupe on se ram\`ene sans perdre en g\'en\'eralit\'e au cas o\`u $G$ 
et $\Gamma$ v\'erifient la propri\'et\'e :
\begin{equation}
\label{eq:H2}
\tag{H2}
\text{\em L'intersection de $\Gamma$ et du centre de $G$ est 
triviale.}
\end{equation}
Notons que la propri\'et\'e~\eqref{eq:H1} reste v\'erifi\'ee apr\`es cette deuxi\`eme 
r\'eduction.

La propri\'et\'e~\eqref{eq:H2} signifie simplement que l'action par translation de $G$ 
sur $X$ est libre, c'est \`a dire que $G$ peut \^etre consid\'er\'e comme un 
groupe de transformations de $X$. Elle entra\^\i ne imm\'ediatement que 
 $\Gamma$ nilpotent d'ordre $k-1$. En particulier, si $k=2$ alors 
 $\Gamma$ est ab\'elien.

 Des rappels pr\'ec\'edents on d\'eduit 
facilement :
\begin{proposition}\label{prop:nilprojectif}
Soient $(X,T)$ un syst\`eme minimal et $k\geq 1$ un entier.
 Alors la famille des facteurs de $X$ 
qui sont isomorphes \`a des nilsyst\`emes d'ordre $k$ est projective.
\end{proposition}
\begin{proof}
Soient $p_1\colon X\to X_1$ et $p_2\colon X\to X_2$ deux facteurs 
de 
$X$ et supposons que $X_1$ et $X_2$ sont des nilsyst\`{e}mes d'ordre $k$.

Soit $Y\subset X_1\times X_2$ l'image de l'application $p_1\times 
p_2\colon X\to X_1\times X_2$.
Alors $Y$ est un ferm\'e de $X_1\times X_2$, invariant par 
$T_1\times T_2$. D'apr\`es les rappels pr\'ec\'edents, $(Y,T_1\times T_2)$ 
est un nilsyst\`eme d'ordre $k$.

De plus, $(Y,T_1\times T_2)$ est un facteur de $(X,T)$ par 
l'application facteur $p_1\times p_2$ et $X_1,X_2$  et les 
applications $p_1$, $p_2$ 
sont les compos\'ees de $p_1\times p_2$ par les projections naturelles 
qui sont des applications facteur.
 Ainsi, $Y$ est un facteur de $X$ \guillemotleft~au 
dessus~\guillemotright\ de $X_1$ et $X_2$, et il est imm\'ediat que c'est le 
plus petit possible.
\end{proof}

\subsection{La relation $\RPD$ dans les nilsyst\`emes d'ordre deux}
Nous nous restreignons d\'esormais aux nilsyst\`emes d'ordre deux.
On rappelle que le groupe $G$ est nilpotent d'ordre $2$ si et
seulement si
son groupe des commutateurs   $G_2$  est inclus dans 
son centre $\CZ(G)$. Les hypoth\`eses~(H1) et~(H2) entra\^\i nent 
 que le groupe de Lie ab\'elien $G_2$ est compact et connexe, c'est donc 
 un tore de dimension finie.
\begin{proposition}
\label{prop:nilsystemes}
Soit $(X,T)$ un nilsyst\`eme minimal d'ordre deux. Alors la relation de 
bi-proximalit\'e r\'egionale sur $X$ est l'identit\'e. De plus,
$(\CP,\CQ)$ est une structure de parall\'el\'epip\`edes forte sur $X$.
\end{proposition}

\begin{proof}
Soit $(X,T)$ un nilsyst\`eme minimal d'ordre deux. 
Soient $G$, $\Gamma$, $t$ et $T$ comme dans la 
d\'efinition~\ref{def:nilsysteme} et supposons que les 
propri\'et\'es~\eqref{eq:H1} et~\eqref{eq:H2} sont v\'erifi\'ees.

Soit $(\CP_X,\CQ_X)$ la structure forte de parall\'el\'epip\`edes  
sur $X$  construite dans la section~3.7 de~\cite{HK} en prenant $F=G_2$.
Pour montrer la proposition, il suffit de v\'erifier que
 $\CP_X$ et $\CQ_X$ sont \'egaux aux 
ensembles $\CP$ et $\CQ$ de la d\'efinition~\ref{def:CPCQ}.

Nous utilisons librement les 
notations et r\'esultats de~\cite{HK}.
On rappelle que $\CP_X$ est l'image de 
$G\type{2,1}\Gamma\type 2$
dans $X^4=G\type 2/\Gamma\type 2$. et que 
$$
 \CP_X=\{(x_0,x_1,x_2,x_3)\in X^4\;;\; 
 \pi'(x_0)\pi'(x_1)\inv\pi'(x_2)\inv\pi'(x_3)= 1\}
$$
o\`u $\pi'$ est la surjection naturelle de $X=G/\Gamma$ sur le groupe 
$G/G_2\Gamma$. D'apr\`es le th\'eor\`eme~\ref{th:rappelsnil1}, 
$G/G_2\Gamma$ est le facteur \'equicontinu maximal de $X$ et $\pi'$ 
co•ncide avec l'application facteur $\pi$. D'apr\`es le 
th\'eor\`eme~\ref{th:para-distal}, $\CP_X=\CP$.

D'autre part,  on rappelle que $\Gamma\type{3,1}=\Gamma\type 3\cap 
G\type{3,1}$ et que $\CQ_X$ est par d\'efinition l'image 
de $G\type{3,1}/\Gamma\type {3,1}$ dans 
$X^8=G\type 3/\Gamma\type 3$ par l'inclusion naturelle.
On v\'erifie facilement que 
$\Gamma\type{3,1}$ est cocompact dans $G\type{3,1}$ et on en d\'eduit 
que  $\CQ_X$ est ferm\'e dans $X^8$.

Pour tout $x\in X$ et tous $m,n,p\in \Z$, 
$(x,T^mx,T^nx,\dots,T^{m+n+p}x)$ appartient \`a $\CQ_X$ par d\'efinition, 
donc $\CQ\subset\CQ_X$. Montrons l'inclusion oppos\'ee.
Soit $\ubx=(x_0,\dots,x_6,x_7)\in \CQ_X$. 
Alors  $(x_0,x_1,x_2,x_3)$, $(x_0,x_1,x_4,x_5)$ et
$(x_0,x_2,x_4,x_6)$ appartiennent \`a $\CP_X$ donc \`a $\CP$.
D'apr\`es la proposition~\ref{prop:EFFP}, il existe $x'_7\in X$ tel que 
$\ubx'=(x_0,\dots,x_6,x'_7)$ appartienne \`a $\CQ$. On a donc
$\ubx'\in\CQ_X$, et comme $(\CP_X,\CQ_X)$ est une 
structure forte de parall\'el\'epip\`edes et que $\ubx\in\CQ_X$,
 $x'_7=x_7$ et donc
$\ubx\in\CQ$.
\end{proof}

\begin{remarque} Il est aussi possible de montrer la 
proposition~\ref{prop:nilsystemes} sans 
sans utiliser les  
r\'esultats de la section~\ref{sec:distaux}
(th\'eor\`eme~\ref{th:para-distal} et  
proposition~\ref{prop:EFFP}) mais en d\'ecrivant explicitement 
$\CP$ et $\CQ$. 

Pour d'autres pr\'esentations de $\CP_X$ et $\CQ_X$ on pourra consulter
la section~11 de~\cite{HK2} ou l'appendice~E de~\cite{GT}.
\end{remarque}

\begin{corollaire}
\label{cor:nilsystemes}
Supposons que $(X,T)$ est une limite projective de nilsyst\`{e}mes minimaux
d'ordre deux. Alors la relation de bi-proximalit\'e r\'egionale de $X$ est 
l'\'egalit\'e.
\end{corollaire}

\begin{proof}
Soit $(X_i)$ une famille projective de nilsyst\`{e}mes d'ordre $2$ 
 dont 
la limite projective est $X$. Soient $x,y\in X$ deux points 
bi-r\'egionalement proximaux; D'apr\`es le lemme~\ref{lem:facteur}, pour 
tout $i$ les images de $x$ et $y$ dans $X_i$ sont bi-r\'egionalement 
proximales et donc \'egales d'apr\`es la proposition~\ref{prop:nilsystemes}. 
On a donc $x=y$.
\end{proof}

\subsection{D\'emonstration du th\'eor\`eme~\ref{th:main-distal} \`a partir du 
th\'eor\`eme~\ref{th:main}}
\label{subsec:proof-main-distal}

Soit $(X,T)$ un syst\`eme distal minimal. 
Nous avons vu que $(\CP,\CQ)$ est une structure de parall\'el\'epip\`edes 
sur $X$.
D'apr\`es la section~3.3 de~\cite{HK} , la relation de bi-proximalit\'e 
r\'egionale de $X$ est donc une relation d'\'equivalence sur $X$. 
Soient $Y$ le  quotient de $X$ par cette relation et 
$\phi\colon X\to Y$ l'application quotient.

Comme la relation de bi-proximalit\'e r\'egionale est ferm\'ee, $Y$ peut 
\^etre muni d'une structure d'espace compact telle que $\phi$ soit 
continue et comme cette relation  
 est 
invariante par $T$, cette transformation induit une transformation 
de $Y$, encore not\'ee $T$, telle que $\phi\colon X\to Y$ soit un 
facteur.

D'apr\`es la section~3.3 de~\cite{HK}, l'image de $\CQ(X)$ par 
l'application $\phi\type 2$ est une structure de parall\'el\'epip\`edes 
forte sur $Y$.
Par ailleurs, d'apr\`es le lemme~\ref{lem:facteur},  
cette image est \'egale \`a $\CQ(Y)$ qui est 
donc une structure forte. 

On en d\'eduit que la relation de 
bi-proximalit\'e r\'egionale de $Y$ est l'\'egalit\'e.
D'apr\`es le th\'eor\`eme~\ref{th:main}, $Y$ est une limite projective de 
nilsyst\`{e}mes d'ordre $2$.

Soit maintenant $(W,T)$ un nilsyst\`eme d'ordre $2$ qui est un 
facteur de $X$, l'application facteur \'etant not\'ee $\psi$.
Soient $x,y\in X$ deux points tels que $\phi(x)=\phi(y)$ et montrons 
que $\psi(x)=\psi(y)$. 

Par d\'{e}finition de $\phi$, $x$ et $y$ sont bi-r\'egionalement proximaux 
et il existe $a,b,c\in X$ tels que $(x,y,a,a,b,b,c,c)\in\CQ(X)$. 
On a donc
$(\psi(x),\psi(y),\psi(a),\psi(a),\psi(b),\psi(b),$ $\psi(c),\psi(c))
\in\CQ(W)$ 
donc $\psi(x)$ et $\psi(y)$ sont bi-r\'egionalement proximaux 
dans $W$, donc \'egaux d'apr\`es la proposition~\ref{prop:nilsystemes}. 
Notre affirmation est d\'emontr\'ee.

Ainsi, l'application $\psi$ se factorise \`a travers $\phi$ : tout 
facteur de $X$ qui est un nilsyst\`eme d'ordre $2$ est un facteur de 
$W$, qui est donc le nilfacteur d'ordre $2$ maximal de $X$.
\qed

\section{R\'esultats topologiques}
\label{sec:topo}

Dans cette section nous la d\'ebutons la d\'emonstration du 
th\'eor\`eme~\ref{th:main}.

Nous supposons d\'esormais :

\smallskip
\begin{itemize}
\item[(Hyp.)]
\em Le syst\`eme $(X,T)$ est transitif et la relation 
de bi-proximalit\'e r\'egionale sur $X$ est l'\'egalit\'e.
\end{itemize}
\smallskip
Nous voulons d\'emontrer que 
le syst\`eme est une limite projective de nilsyst\`emes d'ordre deux.
La d\'emonstration occupe le reste de cette section et toute la section 
suivante. 
D'apr\`es la proposition~\ref{prop:distalRP}, le syst\`eme $(X,T)$ est 
 distal et minimal.
D'apr\`es le th\'eor\`eme~\ref{th:struc-distal}, $(\CP,\CQ)$ est une 
structure de parall\'el\'epip\`edes sur $X$ telle qu'elle a \'et\'e d\'efinie et 
\'etudi\'ee dans~\cite{HK}, et cette   structure est forte 
par hypoth\`ese. 

Les structures de parall\'elogrammes et de parall\'el\'epip\`edes ont \'et\'e 
\'etudi\'ees dans cet article dans un cadre  \guillemotleft~abstrait~\guillemotright\ et nous avons 
ici une  structure \guillemotleft~plus riche~\guillemotright\ : 
 D'une part, $X$, $\CP$, 
$\CQ$, \dots\  sont des espaces compacts. D'autre part, les 
parall\'elogrammes et parall\'{e}l\'{e}pip\`{e}des de $X$ proviennent de la dynamique.
Les structures topologiques et dynamiques sont bien s\^ur intimement 
li\'ees mais, pour la commodit\'e de la lecture, nous s\'eparons les deux 
type d'arguments et  commen\c cons par ceux 
de nature purement topologique.

Nous utilisons d\'esormais librement les d\'efinitions, notations et 
r\'{e}sultats de~\cite{HK}.

\subsection{Les groupes $P_s$ et $F$.}
Nous \'etudions les propri\'et\'es topologiques des groupes et 
homomorphismes introduites dans la section~4 de~\cite{HK}. Les 
d\'emonstrations consistent essentiellement en arguments \'el\'ementaires de 
compacit\'e que
nous ne donnons que succinctement. 

Rappelons que la classe d'un parall\'elogramme $\bx$ pour la
 relation d'\'equivalence $\approx$ est not\'ee $\br\bx$ et que 
$P$ est le quotient de $\CP$ par cette relation. 

 On a $\CQ=\{(\bx,\by)\;;\;\bx,\by\in\CP,\ 
\bx\approx\by\}$. Ainsi, le graphe de la relation $\approx$ est le 
ferm\'e $\CQ$.
On peut donc munir $P$ d'une structure d'espace compact telle que la 
surjection naturelle $\CP\to P$ soit continue.

On d\'efinit deux applications $r,q\colon P\to Z$ par
$$
\text{pour tout }\bx\in\CP,\ 
 r(\br\bx)=\pi(x_1)\pi(x_0)\inv\text{ et 
 }q(\br\bx)=\pi(x_2)\pi(x_0)\inv\ .
$$
Ces deux applications sont continues. 
En particulier, pour tout $s\in Z$, le sous-ensemble 
$P_s=r\inv(\{s\})$ de $P$ est ferm\'e, et la restriction $q_s$ de $q$ \`a 
$P_s$ est continue. 

Dans la section~4.1 de~\cite{HK} on a muni chaque ensemble $P_s$ d'une structure de groupe 
ab\'elien. Montrons :

\begin{lemme}
\label{lem:m-continue}
Pour chaque $s\in Z$ la multiplication de $P_s$ est continue. Plus 
pr\'ecis\'ement, soit 
$$
 P\times_r P=\bigl\{(\alpha,\beta)\in P\times P\;;\;r(\alpha)=r(\beta)\bigr\}
=\bigcup_{s\in Z}P_s\times P_s\ .
$$
Soit  $m\colon  P\times_r P\to P$ l'application dont, pour tout $s\in Z$,
la restriction \`a $P_s\times P_s$ est la multiplication. Alors $m$ est 
continue.
\end{lemme}
\begin{proof}
D\'efinissons
$$
 K=\bigl\{(x_0,x_1,x_2,x_3,x_4,x_5)\in X^6\;;\; 
 (x_0,x_1,x_2,x_3)\in\CP\text{ et }(x_2,x_3,x_4,x_5)\in\CP\bigr\}\ .
$$
Ainsi, $K$ est une partie ferm\'ee de $X^6$ et $P\times_r P$ est l'image de $K$ 
par l'application 
$$
(x_0,x_1,x_1,x_3,x_4,x_5)\mapsto \bigl(\br{x_0,x_1,x_2,x_3},
\br{x_2,x_3,x_4,x_5}\bigr)\ .
$$
L'application $m$ de $P\times_r P$ dans $P$ est caract\'eris\'ee par :
$$
 \text{pour tout }(x_0,x_1,x_1,x_3,x_4,x_5)\in R,\ 
m\bigl(\br{x_0,x_1,x_2,x_3},
\br{x_2,x_3,x_4,x_5}\bigr)=\br{x_0,x_1,x_4,x_5}\ .
$$
Le graphe de cette application est ferm\'e et elle est donc continue.
\end{proof}
Ainsi, pour tout $s\in Z$, $P_s$ est un groupe ab\'elien compact et 
$q_s\colon P_s\to Z$ est un homomorphisme continu de  groupes.
Le noyau $F_s$ de cet homomorphisme est donc un groupe compact. 

Notons $\CP_{1,1}$ l'ensemble parall\'{e}logrammes verticaux,  c'est \`a 
dire dont les quatre sommets ont la m\^eme projection sur $Z$.
On rappelle que $F$ est form\'e des classes d'\'equivalence des 
parall\'{e}logrammes verticaux. Ainsi, $F=F_1$ est un groupe compact.

\begin{lemme}
L'action de $F$ sur $X$ est continue.
\end{lemme}
\begin{proof}
 D\'efinissons un sous-ensemble $K$ de 
$X\times X\times \CP_{1,1}$ par :
$$
 K=\bigl\{(x,y,\bz)\;;\; x,y\in X,\ \bz\in\CP_{1,1},\ 
 \pi(x)=\pi(y)\text{ et }
(x,x,x,y,\bz\in\CQ\bigr\}\ .
$$
Alors $K$ est ferm\'e dans $X\times X\times\CP_{1,1}$ donc compact. 
Remarquons que si $(x,y,\bz\bigr)\in K$, si $\bz'\in\CP_{1,1}$ 
et si $\br{\bz'}=\br{\bz}$ alors $(x,y,\bz'\bigr)\in K$. Ainsi, 
$K$ induit un sous-ensemble compact $R$ de $X\times X\times F$. 

Par d\'efinition de l'action de $F$ sur $X$, pour $x\in X$ et $u\in F$, 
$u\cdot x$ est l'unique \'el\'ement de $X$ tel que 
$\bigl(x,u\cdot x,u)$ appartienne \`a $R$. Cette action est donc 
continue.
\end{proof}

Introduisons quelques notations suppl\'ementaires.
\begin{notation}On note
$$
 L=q\inv(\{1\})=\bigcup_{s\in Z} F_s\subset P
$$
et $j\colon L\to F$ l'application donn\'ee par 
$$
j(\alpha)=j_s(\alpha)\text{ si }\alpha\in F_s
$$
Pour chaque $s\in Z$ on note $i_s\colon F\to F_s$ l'application 
r\'eciproque de $j_s$ et  $i\colon Z\times F\to L$ d\'esigne
l'application d\'efinie par $i(s,\xi)=i_s(\xi)$.
\end{notation}

\begin{lemme}
Pour tout $s\in Z$, l'isomorphisme $j_s\colon F_s\to F$ est continu et 
plus pr\'{e}cis\'{e}ment l'application $j\colon L\to F$
est continue.
\end{lemme}
\begin{proof}
Soit $\CP_{\cdot,1}$ l'ensemble des parall\'elogrammes 
$\bx=(x_0,x_1,x_2,x_3)$ v\'erifiant $\pi(x_2)=\pi(x_0)$ et donc aussi
$\pi(x_3)=\pi(x_1)$. 
Ainsi, $L$ est l'image de $\CP_{\cdot,1}$ dans le quotient $P$ de $\CP$.

Tout 
parall\'elogramme appartenant \`a $\CP_{\cdot,1}$ peut s'\'ecrire d'une 
unique mani\`ere sous la forme $(a,b,u\cdot a,v\cdot b)$ avec $a,b\in X$ 
et $u,v\in F$ et on d\'efinit une application $j'\colon \CP_{\cdot,1}\to 
F$ par 
$j'(a,b,u\cdot a,v\cdot b)=vu\inv$.
Cette application est continue. De plus, $j$ est l'application induite 
par $j'$ par passage au quotient, et $j$ est donc continue.
\end{proof}
Remarquons que l'application $i$ est la r\'eciproque de l'application 
$r\times j\colon L\to Z\times F$. On a donc:
\begin{corollaire}\label{cor:icontinue}
L'application $i\colon Z\times F\to L$ est continue.
\end{corollaire}

\subsection{Le groupe de structure topologique}
\begin{definition}
Le \emph{groupe de structure topologique} de $(X,T)$, not\'e $\CGtop$ 
ou $\CGtop(X)$, est l'ensemble des hom\'eomorphismes de $X$ appartenant 
au groupe de structure $\CG$ de $X$.
\end{definition}
Ainsi, $\CGtop$ est le groupe form\'e des hom\'eomorphismes $x\mapsto 
g\cdot x$ de $X$ tels que pour tout $\bx\in\CP$ on ait $g\type 
2\cdot\bx\in\CP$ et $g\type 2\cdot\bx\equiv\bx$. On rappelle que
$$
\text{pour }\bx=(x_0,x_1,x_2,x_3)\in\CP,\ g\type 2\cdot\bx=
(g\cdot x_0,g\cdot x_1,g\cdot x_2,g\cdot x_3)\ .
$$
On munit $\CGtop$ de la topologie de la convergence uniforme sur $X$. 
Ainsi, $\CGtop$ est un groupe polonais. Comme $\CGtop$ est un 
sous-groupe de $\CG$, il est nilpotent d'ordre $2$ (\cite{HK}, 
proposition~12).

\begin{proposition}\label{prop:transitive2}
Le groupe $\CGtop$ agit transitivement sur $X$ si et seulement si, 
pour tout $s\in Z$ la suite exacte
\begin{equation}
\label{eq:exacttopo}
0\longrightarrow F\overset{i_s}{\longrightarrow} P_s
\overset{q_s}{\longrightarrow} B
\longrightarrow 0
\end{equation}
se scinde contin\^ument, c'est \`a dire qu'il existe
un homomorphisme de groupes continu 
$\kappa_s\colon P_s\to F$ tel que $\kappa_s\circ j_s$ soit l'identit\'e 
de $F$.

En particulier, si $F$ est un tore de dimension finie alors $\CGtop$ 
agit transitivement sur $X$.
\end{proposition}

\begin{proof} Il suffit de v\'erifier que l'\'el\'ement de $\CG$ construit 
dans la preuve du th\'eor\`eme~2 de~\cite{HK} est un hom\'eomorphisme de $X$.
\end{proof}

\begin{proposition}
\label{prop:transitive}
Supposons que $\CGtop$ agit 
transitivement sur $X$.
Soit $\Gammatop$ le stabilisateur d'un point $e$ de $X$ et identifions 
$X$ avec $\CGtop/\Gammatop$ de la fa\c con naturelle. 

Alors 
$\CGtop$ est localement compact, $\Gammatop$ est un sous-groupe discret de 
$G$, et l'identification $X=\CGtop/\Gammatop$ est un hom\'eomorphisme.
\end{proposition}
Avant la d\'emonstration de cette proposition, introduisons une notation utilis\'ee dans toute la suite.
\begin{notation}
Si $H,K$ sont deux groupes topologiques m\'etrisables, 
$\hom_c(H,K)$ d\'esigne 
l'ensemble des homomorphismes continus de $H$ dans $K$, muni de la 
topologie de la convergence uniforme sur tous les compacts de $H$.
\end{notation}

\begin{proof} Pour all\'eger la typographie, dans la d\'emonstration nous 
\'ecrivons $G$ et $\Gamma$ au lieu de $\CGtop$ et $\Gammatop$.

\subsubsection{Remarques pr\'eliminaires.}

 Comme $G$ est un groupe transitif de transformations de $X$
et que $\Gamma$ est le stabilisateur d'un point $e$, l'intersection 
de ce groupe et du centre $\CZ(G)$ de $G$ est triviale. D'apr\`es la 
proposition~15 de~\cite{HK}, $F$ est inclus dans $\CZ(G)$ donc 
$\Gamma\cap F=\{1\}$. Dans la section~5.3 de~\cite{HK} on montre que 
le groupe $G_2$ des commutateurs de $G$ est inclus dans $F$. On a donc
$\Gamma\cap G_2=\{1\}$ et $\Gamma$ est ab\'elien.

Par d\'efinition, $\Gamma$ est un sous-groupe ferm\'e de $G$. 
Comme $F$  agit contin\^{u}ment sur $X$, l'inclusion de $F$ dans $G$ est continue, 
et $F$ est un sous-groupe compact de $G$.

\subsubsection{Premi\`ere partie.}
On rappelle encore (\cite{HK}, section~5.3) qu'il existe un 
homomorphisme de groupes $p\colon G\to Z$, v\'erifiant
$$
\text{pour tout $g\in G$ et tout $x\in X$, } p(g)\cdot\pi(x)=\pi(g\cdot x)\ .
$$
Par d\'efinition de l'identification $X=G/\Gamma$, le point $e$ de $X$ 
est l'image dans $X$ de l'\'el\'ement unit\'e $1$ de $G$,  donc $\pi(e)$ 
 est l'\'el\'ement unit\'e $1$ de $Z$. Ainsi, 
$$
 p(g)=\pi(g\cdot e)\text{ pour tout }g\in G\ .
$$

Le graphe de l'application  $p$ est ferm\'e et cet homomorphisme est donc 
continu. 
Son noyau $F\Gamma$ est donc ferm\'e dans $G$.
 Comme $G$ et $Z$ sont des groupes polonais et 
que $p$ est surjectif, cet 
homomorphisme est une application ouverte  et la topologie de $Z$ 
co•ncide avec la topologie quotient de $G/F\Gamma$ (voir~\cite{BK}, 
chapitre 1).

La compacit\'e de $F$ entra\^\i ne alors facilement que l'application 
$g\mapsto g\cdot e$ est ouverte de $G$ dans $X$. 
Ainsi, l'identification $X=G/\Gamma$ est un hom\'eomorphisme et $\Gamma$ 
est un sous-groupe cocompact de $G$.

\subsubsection{Deuxi\`eme partie.} Nous montrons maintenant que 
$\Gamma$ est un sous-groupe  discret de $G$.

\noindent a) Soit $\gamma\in\Gamma$. Comme $G$ est nilpotent 
d'ordre $2$, 
l'application 
$$f_\gamma\colon g\mapsto [\gamma,g]$$
 est un homomorphisme de groupes 
de $G$ dans $G_2\subset F$. 

Comme l'application $(h,g)\mapsto[h,g]\colon G\times G\to G_2$ est 
continue, pour chaque $\gamma\in\Gamma$ l'homomorphisme $f_\gamma$ 
est continu et
 l'application $\gamma\mapsto 
f_\gamma$ est continue de $\Gamma$ dans $\hom_c(G,F)$.

Soit $\gamma\in\Gamma$. Comme $F$ est contenu dans le centre de $G$, 
l'homomorphisme $f_\gamma$ est trivial sur $F$. Comme $\Gamma$ est 
ab\'elien, la restriction de l'homomorphisme $f_\gamma$ \`a $\Gamma$ est 
triviale. Ainsi, $f_\gamma$ induit un homomorphisme de groupes
 $\chi_\gamma$ de $G/\Gamma F=Z$ dans $F$, et cet homomorphisme est 
continu. 

On v\'erifie facilement que $\chi_\gamma\neq 1$ pour tout 
$\gamma\in\Gamma$ diff\'erent de $1$. En effet, dans le cas contraire on 
aurait $f_\gamma=1$ et $\gamma$ appartiendrait au centre de $G$ 
ce qui 
contredit les remarques pr\'eliminaires.

\noindent b)
Supposons que $\Gamma$ n'est pas discret. 

Il existe alors une suite $(\gamma_n)$ 
dans $\Gamma$, convergeant vers $1$ dans $\Gamma$ et avec 
$\gamma_n\neq 1$ pour tout $n$.

On a donc $\chi_{\gamma_n}\neq 1$ pour tout $n$.
Comme $Z$ et $F$ sont des groupes ab\'eliens compacts, le groupe 
$\hom_c(Z,F)$ est discret et la suite $(\chi_{\gamma_n})$ ne peut pas 
converger vers $1$ dans $\hom_c(Z,F)$, c'est \`a dire uniform\'ement sur 
$Z$. 

Il existe donc  une suite $(z_n)$ dans $Z$ telle 
que la suite $(\chi_{\gamma_n}(z_n))$  converge dans $F$ vers un 
point $u\neq 1$.
 En repla\c cant la suite $(\gamma_n)$ par une sous-suite on se 
am\`ene au cas o\`u la suite $(z_n)$ converge dans $Z$ vers un certain point 
$z$.
Comme $p$ est une application ouverte, il existe une suite 
 $(g_n)$ 
dans $G$ avec $p(g_n)=z_n$ pour tout $n$ et qui converge vers un 
point $g$ de $G$ avec $p(g)=z$. 

Le sous-ensemble $K=\{g_n\;;\;n\geq 1\}\cup\{g\}$ est un compact de $G$. 
Comme l'application $\gamma\mapsto f_\gamma$ est continue \`a valeurs 
dans le 
groupe $\hom_c(G,F)$  muni de la topologie de la convergence 
uniforme sur tout compact, la suite $(f_{\gamma_n}(\cdot))$ converge 
uniform\'ement vers $1$ sur $K$ et la suite $(f_{\gamma_n}(g_n))$ 
converge donc vers $1$. 

Comme pour tout $n$ on a $f_{\gamma_n}(g_n)=\chi_{\gamma_n}(z_n)$ on 
obtient la contradiction recherch\'ee.

\subsubsection{Troisi\`eme partie.} Nous montrons maintenant  que $G$ est 
localement compact.

 Comme $\Gamma$ est un sous-groupe discret de $G$, il existe 
un voisinage ferm\'e sym\'etrique $K$ de $1$ dans $G$  avec $(K\cdot 
K\cdot K)\cap\Gamma=\{1\}$. V\'erifions que $K$ est compact.
Soit $(g_n)$ une suite dans $K$ et montrons que cette suite admet une 
sous-suite convergente. 

Comme $X$ est compact, quitte \`a remplacer la 
suite donn\'ee par un sous-suite, on peut supposer que la suite 
$(g_n\cdot e)$ converge dans $X$. Comme l'application $g\mapsto g\cdot 
e$ est ouverte, il existe une suite $(h_n)$ dans $G$, convergente et 
v\'erifiant $h_n\cdot e=g_n\cdot e$ pour tout $n$. 
Pour tout $n$ posons $\gamma_n=g_n\inv h_n$. On a $\gamma_n\in\Gamma$.

Pour tous $m,n$ assez grands, $h_m h_n\inv\in K$ donc 
$\gamma_m\gamma_n\inv\in K\cdot K\cdot K$ et donc $\gamma_m=\gamma_n$. 
Ainsi, la suite $(\gamma_n)$ est constante \`a partir d'un certain rang 
et la suite $(g_n)$ converge.
\end{proof}

\section{D\'emonstration du th\'eor\`eme~\ref{th:main}}
\label{sec:dyna}
Dans cette section nous terminons la d\'emonstration du
th\'eor\`eme~\ref{th:main}. Nous supposons toujours que
\smallskip
\begin{itemize}
\item[(Hyp.)]
\em La relation 
de bi-proximalit\'e r\'egionale sur $X$ est l'\'egalit\'e.
\end{itemize}
\smallskip

Nous commen\c cons par une remarque pr\'eliminaire.

\subsection{Passage au quotient.}
\label{subsec:quotient}
Soit $F_0$ un sous-groupe ferm\'e de $F$ et 
 $Y$ l'espace quotient de $X$ sous l'action de $F_0$.

 Comme la transformation $T$ appartient au 
groupe de structure de $X$ et que $F$ est contenu dans le centre de ce 
groupe, l'action de $F$ sur $X$ commute avec $T$. Ainsi $T$ induit un 
hom\'eomorphisme, encore not\'e $T$,  de $Y$ et l'application quotient 
$\phi\colon X\to Y$ est un facteur. 
\begin{lemme}
\label{lem:quotient}
Soient $F_0$ et $Y$ comme plus haut. Alors :
\begin{enumerate}
\item\label{it:quotient1}
la relation de bi-proximalit\'e r\'egionale sur $Y$ est 
l'\'egalit\'e et
\item\label{it:quotient2}
le groupe des fibres de $Y$ est isomorphe \`a 
$F/F_0$.
\end{enumerate}
\end{lemme}
\begin{proof} 
\ref{it:quotient1}
 Soient $y,y'\in Y$ avec $(y,y')\in\RPD(Y)$. Alors d'apr\`es 
le lemme~\ref{lem:RPRQ} il
existe $a,b,c\in Y$ avec $(y,y',a,a,b,b,c,c)\in\CQ(Y)$.

D'apr\`es le lemme~\ref{lem:facteur},
$\CQ(Y)$ est l'image de $\CQ(X)$ par $\phi\type 3$ et il 
existe donc $\ubx=(x_0,\dots,x_7)\in\CQ(X)$ avec
$$
 \phi(x_0)=y,\ \phi(x_1)=y',\ \phi(x_2)=\phi(x_3)=a,\ 
 \phi(x_4)=\phi(x_5)=b\text{ et }\phi(x_6)=\phi(x_7)=c\ .
$$
Par d\'efinition de $\phi$, il existe $u,v,w\in F_0$ avec
$$
 x_3=u\cdot x_2,\ x_5=v\cdot x_4\text{ et }x_7=w\cdot x_6\ .
$$
D'apr\`es la proposition~12 de~\cite{HK} on a donc $x_1=uvw\inv\cdot x_1$ et donc par d\'efinition 
de $\phi$, $\phi(x_1)=\phi(x_0)$ c'est \`a dire $y'=y$, ce qui 
d\'emontre~\ref{it:quotient1}.

\medskip\noindent\ref{it:quotient2}
L'image par l'application $\phi\type 2\colon X\type 2\to Y\type 2$ de 
l'ensemble des parall\'{e}logrammes verticaux de $X$ est 
l'ensemble des parall\'elogrammes verticaux de $Y$. Cette application 
respecte l'\'equivalence des parall\'elogrammes et donc elle induit une 
application surjective $\phi_*$ du groupe des fibres $F$ de $X$ sur 
le groupe des fibres de $Y$. $\phi_*$ est continue puisque 
$\phi\type 2$ est continue et on v\'erifie imm\'ediatement que $\phi_*$ 
 est un homomorphisme de groupes.

Nous remarquons aussi que $\phi_*$ commute avec les actions des 
groupes de fibres sur $X$ et $Y$ :
pour tout $x\in X$ et tout $u\in F$, $\phi(u\cdot 
x)=\phi_*(u)\cdot\phi(x)$.
Le noyau de $\phi_*$ est donc l'ensemble de $u\in F$ tels que 
$\phi(u\cdot x)=\phi(x)$ pour tout $x\in X$, c'est \`a dire $F_0$. 
La propri\'et\'e~\ref{it:quotient2} est d\'emontr\'ee.
\end{proof}

\subsection{Le cas o\`u $F$ est un groupe de Lie}
Rappelons que le groupe ab\'elien compact $F$ est un groupe de Lie si et 
seulement si il peut se repr\'esenter comme la somme directe d'un tore de 
dimension finie  et d'un groupe fini (chacun de ces deux groupes 
pouvant \^etre trivial).

\begin{proposition}
\label{prop:splite-local2}
Supposons que le groupe $F$ est un  groupe de Lie. 
Alors le groupe de structure topologique $\CGtop$ agit transitivement 
sur $X$.
\end{proposition}

\begin{proof}

D'apr\`es la proposition~\ref{prop:transitive2} sommes ramen\'es \`a montrer 
que
pour tout $s\in Z$  la  
suite exacte 
\begin{equation*}
\tag{\ref{eq:exacttopo}}
0\longrightarrow F\overset{i_s}{\longrightarrow} P_s
\overset{q_s}{\longrightarrow} B
\longrightarrow 0
\end{equation*}
est contin\^ument scind\'ee, c'est \`a dire qu'il existe un 
homomorphisme continu $\kappa_s\colon P_s\to F$ v\'erifiant
$$
 i_s\circ\kappa_s=\id_F\ .
$$
La d\'emonstration comprend trois parties.
Nous commen\c cons par v\'erifier que cette suite exacte 
est contin\^ument scind\'ee 
pour $s=1$, puis pour tous les $s$ dans un voisinage $W$ de $1$ dans 
$Z$ et   enfin nous concluons par un argument de minimalit\'e.

\subsubsection{Premi\`ere partie : le cas o\`u $s=1$.}\strut

Soient $x,y,x',y'\in X$ avec $\pi(y)\pi(x)\inv=\pi(y')\pi(x')\inv$. 
Alors $(x,y,x',y')\in\CP$ et d'apr\`es le lemme~9 de~\cite{HK} les 
classes $\br{x,x,y,y}$ et $\br{x',x',y',y'}$ des parall\'elogramme 
$(x,x,y,y)$ et $(x',x',y',y')$ sont identiques. 
 On peut donc d\'efinir une application 
$\rho\colon Z\to P$ par
$$
 \rho(t)=\br{x,x,y,y}\text{ pour tout }(x,y)\in X\times X\text{ tel 
 que } \pi(y)\pi(x)\inv=t\ .
$$
Par construction, l'image de l'application $\rho$ est contenue dans 
$P_1$ et $q_1\circ \rho\colon Z\to Z$ est l'application 
identit\'e. 
L'application $\rho$ est un homomorphisme de  groupes par d\'efinition de 
la multiplication dans $P_1$. Son graphe est l'image de $X\times X$ 
par l'application continue
$$
 (x,y)\mapsto\bigl( \pi(y)\pi(x)\inv, \br{x,x,y,y}\bigr)\colon 
 X\times X \to Z\times P_1
$$
et est donc ferm\'e dans $Z\times P_1$. Cette application est donc 
continue et le groupe $P_1$ est donc bien la somme directe topologique 
 $F$ et $Z$. La suite exacte~\eqref{eq:exacttopo} est donc contin\^ument 
 scind\'ee pour $s=1$ et il existe un homomorphisme continu 
 $\kappa_1\colon P_1\circ F$ avec $i_1\circ\kappa_1=\id_F$.

\subsubsection{Deuxi\`eme partie : construction d'un voisinage de $1$ dans $Z$.}\strut

Nous \'ecrivons $F=\T^d\times K$, o\`u $d\geq 0$ est un entier et $K$ un 
groupe fini. Pour tout $s\in Z$, l'homomorphisme continu injectif 
$i_s\colon T^d\times K\to P_s$ s'\'ecrit alors
$$
 i_s(u,\xi)=\phi_s(u)\cdot \psi_s(\xi)
$$
o\`u $\phi_s\colon\T^d\to P_s$ et $\psi_s\colon K\to P_s$ sont 
des homomorphismes continus injectifs avec
\begin{equation}
\label{eq:alphasbetas}
 \phi_s(\T^d)\cap\psi_s(K)=\{1\}\ .
\end{equation}

\medskip\noindent\emph{a)} Pour tout $s\in Z$ nous construisons un 
homomorphisme continu $\tau_s\colon P_s\to\T^d$ tel que
\begin{equation}
\label{eq:taus}
\tau_s\circ\phi_s\text{ est l'identit\'e de }\T^d\text{ et }
\tau_s\circ\psi_s\colon K\to\T^d\text{ est l'homomorphisme trivial.}
\end{equation}
Soit $s\in Z$. En composant $\phi_s$ avec 
l'application quotient  on obtient un homomorphisme 
continu $\tilde\phi_s\colon \T^d\to P_s/\psi_s(K)$, qui est 
injectif d'apr\`es~\eqref{eq:alphasbetas}. Comme $\T^d$ est un groupe 
injectif, il existe un homomorphisme continu
 $\tilde\tau_s\colon  P_s/\psi_s(K)\to \T^d$ tel que
 $\tilde\tau_s\circ\tilde\phi_s$ soit l'identit\'e de $\T^d$. En 
 composant avec l'application quotient on obtient un homomorphisme 
 continu $\tau_s\colon P_s\to\T^d$ v\'erifiant~\eqref{eq:taus}.

\medskip\noindent\emph{b)}
Nous construisons maintenant un voisinage $W$ de $1$ dans $Z$ et pour tout $s\in 
W$ un homomorphisme continu $\sigma_s\colon P_s\to K$ v\'erifiant
\begin{equation}
\label{eq:sigmas}
\sigma_s\circ\psi_s=\id_K\ .
\end{equation}

En composant $\kappa_1$ avec la projection 
$F=\T^d\times K\to K$ on obtient un homomorphisme continu 
$\sigma_1\colon P_1\to K$ 
v\'erifiant la propri\'et\'e~\eqref{eq:sigmas} pour $s=1$.

Comme $K$ est fini il existe un voisinage ouvert $U$ de $1$ dans $Z$ 
et une application continue $\sigma\colon r\inv(U)\to K$ qui co•ncide 
avec $\sigma_1$ sur $P_1$. Pour $s\in U$ on note $\sigma_s$ la 
restriction de $\sigma$ \`a $P_s$.

Pour chaque sous-ensemble $Y$ de $Z$ notons 
$$
 R_Y=\bigl\{(\alpha,\beta)\in r\inv(Y)\times r\inv(Y)\colon 
 r(\alpha)=r(\beta\bigr\}=\bigcup_{s\in Y}P_s\times P_s\ .
$$
Ainsi, $R_Y$ est inclus dans l'ensemble $R_Z=P\times_r P$ 
consid\'er\'e dans le lemme~\ref{lem:m-continue}.

Soit $m\colon P\times_r P\to P$ l'application  introduite 
dans ce lemme et d\'efinissons $f\colon R_U\to K$ 
par
$$
 f(\alpha,\beta)=
 \sigma(m(\alpha,\beta))\cdot\sigma(\alpha)\inv\cdot\sigma(\beta)\inv\ .
$$
Ainsi, pour tout $s\in U$ et tous $\alpha,\beta\in P_s$,
$$
f(\alpha,\beta)=\sigma_s(\alpha\beta)\cdot\sigma_s(\alpha)\inv
\cdot\sigma_s(\beta)\inv\ .
$$
Comme $\sigma_1$ est un homomorphisme, $f$ est \'egale \`a $1$ sur 
$P_1\times P_1$.
La continuit\'e de $\sigma$ et le lemme~\ref{lem:m-continue} entra\^\i nent 
que
$f$ est continue sur $R_U$. 
Comme $K$ est fini,
il existe donc un voisinage $V$ de $1$ dans $Z$ avec $V\subset U$, 
tel que cette application soit \'egale \`a $1$ sur $R_V$. 
Ainsi, pour tout $s\in V$, $\sigma_s$ est un homomorphisme continu 
de $P_s$ 
dans $K$.

Enfin, d'apr\`es le corollaire~\ref{cor:icontinue}, 
l'application $(s,\xi)\mapsto \psi_s(\xi)$ est continue de 
$V\times K$ dans $r\inv(V)$ et donc
l'application $(s,\xi)\mapsto \sigma_s\circ \psi_s(\xi)$ est continue 
de $V\times K$ dans $K$. Comme $K$ est fini et que $\sigma_1\circ 
\psi_1$ est l'identit\'e de $K$ il existe un voisinage $W$ de $1$ dans $Z$, 
contenu dans $V$, tel que~\eqref{eq:sigmas} soit v\'erifi\'e pour tout
$s\in W$.

\medskip\noindent c)
Remarquons que comme $\T^d$ est connexe et que $K$ est fini, 
$\sigma_s\circ \alpha_s\colon \T^d\to K$ est trivial pour tout $s\in W$. 
Pour tout $s\in W$ posons 
$$
\kappa_s=\tau_s\times\sigma_s\colon P_s\to 
\T^d\times K=F\ .
$$
Alors $\kappa_s$ est un homomorphisme continu et 
les relations~\eqref{eq:taus} et~\eqref{eq:sigmas} entra\^\i nent que
$\kappa_s\circ i_s=\id_F$ pour tout $s\in W$.

La suite exacte~\eqref{eq:exacttopo} est  contin\^ument scind\'ee pour 
tout $s$ dans le voisinage $W$ de $1$ dans $Z$.

\subsubsection{Troisi\`eme partie : Invariance par $T$.}\strut

Nous montrons maintenant que l'ensemble des $s\in Z$ pour lesquels la 
suite exacte~\eqref{eq:exacttopo} est contin\^{u}ment scind\'ee 
est invariant par $T$, ce qui 
ach\`evera la d\'emonstration par minimalit\'e.

Soit $s\in Z$. Pour tout $\bx=(x_0,x_1,x_2,x_3)\in\CP_s$ on a 
$T\type 2_1\bx=(x_0,Tx_1,x_2,Tx_3)\in\CP_{Ts}$. 
L'application $T\type 2_1\colon \CP_s\to\CP_{Ts}$
est clairement continue.

De plus, si 
$\bx,\by\in\CP_s$ v\'erifient $\bx\approx\by$ alors $(\bx,\by)\in\CQ$ et 
$(T\type 2_1\bx,T\type 2_1\by)=T\type 3_1(\bx,\by)\in\CQ$ donc $T\type 
2_1\bx\approx T\type 2_1\by$. Ainsi, $T\type 2_1\colon\CP_s\to\CP_s$  
induit une application continue, encore not\'ee $T\type 2_1$, de $P_s$ dans 
lui-m\^eme. La d\'{e}finition de la multiplication dans $P_s$ et $P_{Ts}$ 
entra\^\i ne imm\'ediatement que cette application est un homomorphisme de 
groupes. En rempla\c cant $T$ par son inverse dans ce qui pr\'ec\`ede on 
obtient que cette application est un isomorphisme. 

On v\'erifie imm\'ediatement sur les d\'efinitions que $q_{Ts}\circ T\type 2_1=q_s$, que
$T\type 2_1F_s=F_{Ts}$ et que $j_{Ts}\circ T\type 2_1=j_s$ donc 
$T\type 2_1\circ i_s=i_{Ts}$. On a donc un diagramme commutatif
$$
 \begin{matrix}
0 \longrightarrow & F & \overset{i_s}{\longrightarrow}  & P_s & 
 \overset{q_s}{\longrightarrow} & Z & \longrightarrow & 0\\
& \Big\downarrow{\scriptstyle\id} & & \Big\downarrow
{\scriptstyle T\type 2_1 }& & \Big\downarrow{\scriptstyle\id} & & \\
0 \longrightarrow & F &\overset{i_{Ts} }{\longrightarrow}  & P_{Ts} & 
 \overset{q_{Ts}}{\longrightarrow} & Z & \longrightarrow & 0
\end{matrix}
$$
o\`u les lignes horizontales sont les suites exactes~\eqref{eq:exacttopo} 
aux points $s$ et $Ts$.  

Si $s$ appartient \`a $Z'$ la premi\`ere suite exacte est contin\^ument 
scind\'ee, donc aussi la deuxi\`eme et donc $Ts$ appartient \`a $Z'$. 
Ainsi $Z'$ est invariant par $T$ et contient un voisinage de $1$, il 
est donc \'egal \`a $Z$ par minimalit\'e de la rotation $(Z,T)$,
ce qui ach\`eve la d\'emonstration.
\end{proof}

\subsection{Fin de la d\'emonstration du th\'eor\`eme~\ref{th:main}}\strut

\noindent\emph{a) R\'{e}duction au cas o\`u $F$ est un groupe de 
Lie.}
Comme $F$ est un groupe ab\'elien compact, il peut s'\'ecrire 
 comme limite projective d'une suite de groupes de Lie.
 Cela signifie qu'il existe une suite d\'ecroissante 
$(F_n)$ de sous-groupes 
ferm\'es de $F$, avec $\cap_nF_n=\{1\}$ et telle que $F/F_n$ soit 
un groupe de Lie pour tout $n$.

Soit $X_n$ le quotient de $X$ par l'action de $F_n$. 
Comme le syst\`eme $X$ est la limite 
projective de la suite $(X_n)$ de syst\`emes, il suffit de montrer que 
chacun d'entre eux est limite projective d'une suite de 
nilsyst\`{e}mes d'ordre $2$.

D'apr\`es le lemme~\ref{lem:quotient},
 la relation de 
bi-proximalit\'{e} r\'egionale de $X_n$ est l'\'egalit\'e et le groupe des fibres 
de $X_n$ est isomorphe \`a $F/F_n$.  
Ainsi, sans perdre en g\'en\'eralit\'e, nous pouvons nous 
restreindre au cas o\`u $F$ est un groupe de Lie.

\medskip\noindent\emph{b)}
Comme $F$  est un groupe de Lie, d'apr\`es la 
proposition~\ref{prop:splite-local2} le groupe $\CGtop$ agit 
transitivement sur $X$. 
 D'apr\`es la proposition~\ref{prop:transitive}, $\CGtop$ est localement 
compact, le stabilisateur $\Gamma$  de $e$ est discret et $X$ peut
s'identifier \`a $\CGtop/\Gamma$. 
On rappelle que $\CGtop$ est nilpotent d'ordre $2$. 

Le groupe localement compact $\CGtop$ peut s'\'ecrire comme limite 
projective de groupes de Lie~\cite{Mo}.
Plus pr\'{e}cis\'{e}ment, il existe 
une suite d\'ecroissante $(K_n)$ de sous-groupes compacts de $\CGtop$, 
contenus dans le centre de $\CGtop$, avec $\cap_nK_n=\{1\}$ et telle que 
$\CGtop/K_n$ soit un groupe de Lie pour tout $n$.

Pour tout $n$, $\CGtop/K_n\Gamma$ muni de la transformation induite 
par $T$ est un nilsyst\`eme d'ordre $2$, et $X$ est la limite 
projective de ces syst\`emes.\qed

 \section{Application : une caract\'erisation des nilsuites d'ordre 
deux}
\label{sec:nilsuites}

\subsection{La d\'efinition}\label{subsec:def-nilsuites}
Soit $(u_n\;;\;n\in\Z)$ une suite born\'ee de nombres complexes. On 
rappelle que cette suite est \emph{presque p\'eriodique} si la famille de 
ses translat\'es est relativement compacte pour la topologie de la 
convergence uniforme. Cette propri\'{e}t\'{e} est satisfaite si et seulement 
si il existe une rotation minimale $(X,T)$, un point $e\in X$ et une 
fonction continue $f\colon X\to\C$ avec $u_n=f(T^ne)$ pour tout 
$n\in\Z$.

La notion de \emph{nilsuite} g\'en\'eralise la notion de suite presque 
p\'eriodique. Nous reproduisons la d\'{e}finition de~\cite{BHK}.
\begin{definition}
Soit $k\geq 1$ un entier.

 Une suite $u=(u_n\;;\;n\in\Z)$ de nombres 
complexes est une 
\emph{nilsuite de base d'ordre $k$} s'il existe un nilsyst\`eme
 $(X,T)$ d'ordre $k$, un point $e\in X$ et une 
fonction continue $f$ sur $X$ tels que $u_n=f(T^ne)$ pour tout 
$n\in\Z$.

Une \emph{nilsuite d'ordre $k$} est une limite uniforme de 
nilsuites de base d'ordre $k$. 
\end{definition}

Ainsi, les  nilsuites d'ordre 
$1$ sont les suites presque p\'eriodiques.

Dans la d\'efinition d'une nilsuite de base on peut supposer sans
perdre en g\'en\'eralit\'e  que le nilsyst\`eme $(X,T)$
est minimal : 
il suffit en effet de remplacer $X$ par l'orbite ferm\'ee de 
$e$, qui munie de $T$ est un nilsyst\`eme d'ordre $k$ 
(voir~\cite{Le}),
transitif donc minimal.

\begin{lemme}
\label{lem:nilsuite}
Soit $k\geq 1$ un entier. Pour que la suite $u=(u_n\;;\;n\in\Z)$
 de nombres complexes soit une nilsuite d'ordre $k$ il faut et il 
 suffit qu'il existe un syst\`eme $(Y,S)$ minimal qui soit une limite 
 projective de nilsyst\`{e}mes d'ordre $k$, un point $a\in Y$ et une 
fonction continue $h$ sur $X$ tels que $u_n=h(S^na)$ pour tout 
$n\in\Z$.
\end{lemme}
 
\begin{proof}
Si la suite $u$ provient d'une limite projective de nilsyst\`emes comme 
dans l'\'enonc\'e, alors par approximation de la fonction $h$ 
on obtient imm\'ediatement que cette suite est une nilsuite d'ordre $k$.

Supposons maintenant que la suite $u=(u_n\;;\;n\in\Z)$ est limite 
uniforme des nilsuites de base $u^{(i)}$ o\`u 
$u^{(i)}=(u^{(i)}_n\;;\;n\in\Z)$ pour tout $i\geq 1$.
Pour chaque $i$ choisissons un nilsyst\`eme minimal $(X_i,T_i)$, un 
point $e_i\in X$ et une fonction continue $f_i$ sur $X_i$ avec
$$
 u^{(i)}_n=f_i(T_i^ne_i)\text{ pour tous }i\geq 1\text{ et tous 
 }n\geq 1\ .
$$
Pour tout $i\geq 1$ soient $S_i$ la transformation
$T_1\times\dots\times T_i$ de $X_1\times\dots\times X_i$, 
$a_i=(e_1,\dots,e_i)\in X_1\times\dots\times X_i$ 
et $Y_i$ l'orbite ferm\'ee de $a_i$ pour $S_i$.
Alors $(X_1\times\dots\times X_i, S_i)$ est un nilsyst\`eme d'ordre $k$ 
donc $(Y_i,S_i)$ est un nilsyst\`eme d'ordre $k$ (voir~\cite{Le}), 
transitif donc minimal. 
Notons encore $g_i$ la restriction \`a $Y_i$ de la fonction 
$(x_1,\dots,x_i)\mapsto f_i(x_i)\colon 
X_1\times\dots\times X_i\to\C$. 
On a :
$$
 u^{(i)}_n=g_i(S_i^na_i)\text{ pour tous }i\geq 1
\text{ et tous  }n\geq 1\ .
$$
Pour chaque $i$ soit 
$p_i\colon X_1\times\dots\times X_i\times X_{i+1}
\to  X_1\times\dots\times X_i$ 
la projection naturelle.
On a $p_i(a_{i+1})=a_i$ et $p_i\circ S_{i+1}=S_i\circ p_i$ et par 
densit\'e on en d\'eduit que $p_i(Y_{i+1})=Y_i$. 
Ainsi, $p_i\colon(Y_{i+1},S_{i+1})\to(Y_i,S_i)$ est une application 
facteur.

Les syst\`{e}mes $(Y_i,S_i)$ forment avec les applications $p_i$ un syst\`{e}me 
projectif. La limite projective $(Y,S)$  de ce syst\`eme est un syst\`eme 
minimal. 
Notons  $q_i\colon Y\to Y_i$ l'application naturelle, 
$a$ le point de $Y$ d\'efini par $q_i(a)=a_i$ pour tout $i$ 
et $h_i$ la fonction continue $g_i\circ q_i$ sur $Y$. 
Nous avons 
$$
 u^{(i)}_n=h_i(S^na)\text{ pour tous }i\geq 1
\text{ et tous  }n\geq 1\ .
$$
Quand $i\to+\infty$,  $u^{(i)}_n\to u_n$ uniform\'{e}ment, 
donc la suite de fonctions $(h_i)$ converge uniform\'ement sur l'orbite 
$\{S^na\;;\;n\in\Z)$ de $a$. 
Comme cette orbite est dense dans $Y$, la suite de fonctions 
$(h_i)$ converge uniform\'ement sur $Y$. 
Soit $h$ la limite de cette suite. 
$h$ est une fonction continue sur $Y$ et $u_n=h(S^na)$ pour tout $n$.
\end{proof}

\subsection{Une caract\'erisation des nilsuites d'ordre deux}
\begin{theoreme}
\label{th:nilsuite}
Soit $u=(u_n\;;\;n\in\Z)$ une suite born\'ee de nombres complexes.
\begin{itemize}
\item
Pour que la suite $u$ soit presque p\'eriodique il faut et il 
suffit que 
\begin{align}
 &\text{\emph{pour tout  $\epsilon>0$ 
il existe un entier $M\geq 1$ et $\delta>0$
tels que pour tous}}
\\ 
& k,m,n\in\Z
\text{\emph{ on ait :}}
\notag\\
&\text{\quad\emph{si pour tout }}
i\in[k-M,k+M]\text{ on a }
 |u_{i+m}-u_i|<\delta \text{ \emph{et }}
|u_{i+n}-u_i|<\delta \notag\\
& \text{\quad\emph{alors  on a }}
|u_{k+m+n}-u_k|<\epsilon\ .\notag
\end{align}

\item
Pour que la suite $u$ soit une nilsuite d'ordre deux il faut et il suffit 
que 
\begin{align}
\label{eq:nilsuite}
&\text{\emph{pour tout $\epsilon>0$ il existe un 
entier $M\geq 1$ et $\delta>0$
tels que pour tous}}\\
& k,m,n,p\in\Z \text{\emph{ on 
ait:}} \notag\\
&\text{\quad\emph{si pour tout }} i\in[k-M,k+M]\text{\emph{ on a }}\notag\\
\text{\em(*)} &\qquad |u_{i+m}-u_i|<\delta,\ |u_{i+n}-u_i|<\delta,\ 
 |u_{i+m+n}-u_i|<\delta,\ |u_{i+p}-u_i|<\delta,\ 
\notag\\
 & \qquad\qquad |u_{i+m+p}-u_i|<\delta\text{\emph{ et }} |u_{i+n+p}-u_i|<\delta
\notag\\
& \text{\quad\emph{alors  }}
|u_{k+m+n+p}-u_k|<\epsilon\ .\notag
\end{align}
\end{itemize}
\end{theoreme}
Ainsi, les suites presque p\'eriodiques et les nilsuites d'ordre $2$ 
sont caract\'eris\'ees par des propri\'et\'es de r\'egularit\'e arithm\'etique des 
\emph{temps de retour}. L'uniformit\'e en $k,m,n,p$ 
dans~\eqref{eq:nilsuite} peut se traduire en termes topologiques, 
comme dans la d\'efinition classique d'une suite presque p\'eriodique, 
mais nous ne d\'eveloppons pas ces consid\'erations ici.
\begin{proof}
Nous ne montrons que la deuxi\`eme affirmation, la preuve de la premi\`ere 
\'etant similaire et plus simple.

\noindent \emph{a)} Soient $u=(u_n\;;\;n\in\Z)$ une suite born\'ee et  $D\subset\C$ un disque ferm\'e 
contenant $u_n$ pour tout $n$. 
Muni de la topologie de la convergence simple, $D^\Z$ est un espace 
compact m\'etrisable, par exemple au moyen de la 
distance $d$ d\'efinie par
$$
 \text{pour }x=(x_n\;;\;n\in\Z)\text{ et }y=(y_n\;;\;n\in\Z),\
d(x,y)=\sum_{n\in\Z}2^{-|n|}|x_n-y_n|\ .
$$
Nous munissons $D^\Z$ du \emph{d\'ecalage} $T$ d\'efini par
$$
 \text{pour tout }x\in D^\Z,\ (Tx)_n=x_{n+1}\text{ pour tout }n\in\Z\ .
$$
Ainsi, $T$ est un hom\'eomorphisme de $X$. 

Consid\'erons la suite $u$ comme un \'el\'ement de $D^\Z$ et soit $X$ 
l'orbite ferm\'ee de $u$ sous la transformation $T$. Par 
construction, $(X,T)$ est un syst\`eme transitif.

\medskip\noindent \emph{b)} 
Supposons maintenant que l'hypoth\`ese~\eqref{eq:nilsuite} est 
satisfaite. 

Soient $\epsilon>0$ et $M,\delta$ comme dans~\eqref{eq:nilsuite}.
Soient $N\geq 0$ un entier, $k,m,n,p\in\Z$ et supposons que les 
conditions~(*) sont vraies pour tout    $i\in[k-M-N,k+M+N]$.
Alors, en appliquant l'hypoth\`ese aux entiers $k+j$ avec $|j|\leq N$ 
on obtient que, pour tout $i\in [k-N,k+N]$ on a $|u_{k+m+n+p}-u_k|<\epsilon$.

Par d\'efinition de la distance sur $D^\Z$, cette propri\'et\'e peut s'\'ecrire :
\begin{align}\label{eq:epsdelat2}
&\text{pour tout 
$\epsilon>0$ il existe $\delta>0$ tel que, pour tous $k,m,n,p\in\Z$}\\
&\text{si }
 d(T^{k+m}u,T^ku)<\delta, \ d(T^{k+n}u,T^ku)<\delta,\ 
   d(T^{k+m+n}u,T^ku)<\delta,\notag \\
& \qquad d(T^{k+p}u,T^ku)<\delta, \
  d(T^{k+m+p}u,T^ku)<\delta\text{ et }d(T^{k+n+p}u,T^ku)<\delta
\notag\\
&\text{alors }
 d(T^{k+m+n+p}u,T^ku)<\epsilon\ .\notag
\end{align}
Par densit\'e de $(T^ku\;;\;k\in\Z)$ dans $X$ et par continuit\'e de $T$ 
on obtient  que
\begin{align}
\label{eq:epsdeta2}
&\text{pour tout 
$\epsilon>0$ il existe $\delta>0$ tel que, pour tous $m,n,p\in\Z$ 
et tout $x\in X$,}\\
&\text{si }
 d(T^mx,x)<\delta/2,\ d(T^nx,x)<\delta/2,\ d(T^{m+n}x,x)<\delta/2
\notag\\
&\qquad d(T^px,x)<\delta/2,\ d(T^{m+p}x,x)<\delta/2\text{ et }
d(T^{n+p}x,x)<\delta/2\notag\\
&\text{alors }
 d(T^{m+n+p}x,x)<2\epsilon\ .\notag
\end{align}
Par d\'efinition de $\CQ$ nous obtenons que si $(x,x,\dots,x,y)\in\CQ$ 
alors $x=y$. Le lemme~\ref{lem:RPRQ} entra\^\i ne alors que la relation $\RPDS$ 
est l'\'egalit\'e.  D'apr\`es le th\'eor\`eme~\ref{th:main}, le syst\`eme $(X,T)$ 
est donc une limite projective de nilsyst\`emes d'ordre $2$.

L'application $f$ qui \`a chaque $x=(x_n\;;\;n\in\Z)\in X$ associe 
$x_0$ est une fonction continue, et on a $f(T^nu)=u_n$ pour tout $n$. 
La suite $u=(u_n\;;\;n\in\Z)$ est donc une nilsuite d'ordre $2$.

\medskip\noindent\emph{c)}
Nous montrons maintenant que la propri\'et\'e~\eqref{eq:nilsuite} 
est satisfaite par toute nilsuite d'ordre deux. Comme cette propri\'et\'e 
est stable par limite uniforme nous pouvons nous restreindre au cas 
des nilsuites de base.
 
Soient $(Y,S)$ un nilsyst\`eme minimal d'ordre $2$, $e\in Y$ et $f$ une 
fonction continue sur $Y$ avec $u_n=f(S^ny)$ pour tout $n$. 

Soit $D$ un disque ferm\'e de $\C$ contenant toues les valeur de la 
fonction $f$. D\'efinissons une application $p\colon Y\to D^\Z$ par
$$
 \text{pour tout }y\in Y\text{ et tout }n\in\Z,\ 
 \bigl(p(y)\bigr)_n=f(T^ny)\ .
$$
Cette application v\'erifie clairement $p\circ S=T\circ p$ et $p(e)=u$. 
Par densit\'e, l'image de $p$ est $X$ et  $p$ 
est donc une application facteur de $(Y,S)$ sur $(X,T)$.
Ainsi, $(X,T)$ est un facteur d'un nilsyst\`eme d'ordre $2$ et est donc 
aussi un nilsyst\`eme d'ordre deux, transitif donc minimal.

La relation $\RPDS$ sur $X$ est donc l'\'egalit\'e. On en d\'eduit 
imm\'{e}diatement~\eqref{eq:epsdeta2} puis~\eqref{eq:epsdelat2} et 
enfin~\eqref{eq:nilsuite}.
\end{proof}

\end{document}